\documentclass[12pt]{article}
\usepackage{a4wide}

\usepackage{amssymb,amsfonts,amsmath}
\usepackage{latexsym}
\usepackage[latin1]{inputenc}
\usepackage{comment}
\usepackage[hyperref]{ntheorem}

\usepackage[hyperindex,colorlinks,backref]{hyperref}

\usepackage{amssymb,amsfonts,amsmath}
\usepackage[hyperref]{ntheorem}
\makeatletter

\newtheoremstyle{longplain}
{\item[\hskip\labelsep \theorem@headerfont ##1\
##2\theorem@separator]}
{\item[\hskip\labelsep  \theorem@headerfont ##1\ ##2]
\textbf{(##3)}  \theorem@separator\ \ }

\makeatother

\newtheoremlisttype{list}
{\begin{trivlist}\item}
{\item[\textbf{ ##1 ##2:}]\ ##3\dotfill ##4 }
{\end{trivlist}}

\newtheoremlisttype{tab}
{\begin{tabular*}{\linewidth}{@{}ll@{\extracolsep{\fill}}r@{}}}
{\textbf{##1}\ \textbf{##2}&##3\dotfill&##4\\}
{\end{tabular*}}

\theoremstyle{longplain}
 \newtheorem{thm}{Theorem}[section]
\newtheorem{cor}[thm]{Corollary}
\newtheorem{lem}[thm]{Lemma}
\newtheorem{prop}[thm]{Proposition}
\newtheorem{defn}[thm]{Definition}

\newtheorem{rem}[thm]{Remark}
\newtheorem{con}[thm]{Assumptions}

\numberwithin{equation}{section}

\newcommand{\mymarginpar}[1]{\marginpar{\fbox{
    \begin{minipage}{2cm}
#1
\end{minipage}}}}

\renewcommand{\mymarginpar}[1]{}

\newcommand{\neucom}[2]{
 \marginpar{\fbox{
     \begin{minipage}{1.7cm}
       \textbf{#1} !
     \end{minipage}}}
 \begin{quote}
   \fbox{
     \begin{minipage}{0.8\textwidth} {\color{red}{\it #2}}
     \end{minipage}
   }
 \end{quote}
}

\usepackage{comment}

\renewcommand{\neucom}[2]{}

\usepackage{color}
\definecolor{orange}{rgb}{1.00,0.65,0.00}
\definecolor{purple}{rgb}{0.63,0.13,0.94}

\newcommand{\ra}{{\mathord{\rangle}}}
\newcommand{\la}{{\mathord{\langle}}}
\newcommand{\La}{{\mathord{\Lambda}}}
\newcommand{\g}{{\mathord{\bf g}}}
\newcommand{\K}{{\mathord{\bf K}}}
\newcommand{\h}{{\mathord{\bf h}}}
\newcommand{\e}{{\mathord{\bf e}}}
\newcommand{\pa}{{\mathord{\partial}}}

\newcommand{\setR}{{\mathord{\mathbb R}}}

\begin{document}
\title{On the Well-Posedness of the  Vacuum Einstein's  Equations}
\author{Lavi Karp
\thanks{Research partially supported by DFG: SCHU 808/19-1 and ORT
Braude College's Research Authority}
}

\date{}

\maketitle

\begin{abstract}

The  Cauchy problem of the  vacuum Einstein's  equations aims to find
a semi-metric $\g_{\alpha\beta}$ of a spacetime with
vanishing Ricci
curvature ${\bf R}_{\alpha,\beta} $ and  prescribed initial data.
Under the harmonic gauge condition, the  equations  ${\bf
R}_{\alpha,\beta}=0 $ are transferred  into a system
of quasi-linear wave equations which are called {\it the  reduced
Einstein equations}. The initial data for Einstein's equations are a
proper Riemannian metric $\h_{ab}$ and a second
fundamental form $\K_{ab}$. A necessary condition for the reduced
Einstein equation to satisfy the vacuum equations is that the initial
data satisfy Einstein constraint equations. Hence the
 data  $(\h_{ab},\K_{ab})$ cannot serve as initial data for the
reduced Einstein equations.

 Previous results  in the case of asymptotically flat spacetimes
provide a solution to the constraint equations in one type of Sobolev
spaces, while initial data for the evolution equations  belong to a
different type of Sobolev spaces. The goal of our work is to resolve
this incompatibility and to show that under the harmonic
gauge  the vacuum
Einstein equations are  well-posed in one type of Sobolev
spaces.

\end{abstract}

\maketitle

\section{Introduction}
\label{sec:Introduction}

This paper deals with  well-posedness of the Cauchy problem of
Einstein vacuum equations
\begin{equation}
  \label{eq:1}
  {\bf R}_{\alpha\beta}(\g)=0,\qquad \alpha,\beta=0,1,2,3.
\end{equation}
Here ${\bf R}_{\alpha\beta}(\g)$  denotes the Ricci curvature
tensors of a Lorentzian metric $\g$. The
unknowns are the coefficients $\g_{\alpha\beta}$ of the semi-metric
$\g$.

The initial data consist of the triple $(M,\h_{ab},\K_{ab})$, where
$M$ is a space-like manifold, $\h_{ab}$ is a proper Riemannian metric
on $M$ and $\K_{ab}$ is its second fundamental form (extrinsic
curvature). The
semi-metric $\g_{\alpha\beta}$ takes the following data on $M$:
\begin{equation}
\label{eq:2}
 \left\{\begin{array}{c}
\g_{00}{\mid_{M}}=-1 ,\qquad \g_{0a}{\mid_{M}}=0,
\qquad \g_{ab}{\mid_{M}}=\h_{ab} \\
 -\frac 1 2\partial_0 \g_{ab}{\mid_{M}}=\K_{ab},   \\
\end{array}\qquad a,b=1,2,3.\right.
\end{equation}
The rest of the data $\pa_0 \g_{\alpha0}$ are determined through the
conditions
\begin{equation}
\label{eq:9}
 F^\mu=\g^{\beta\gamma}\Gamma_{\beta\gamma}^\mu=0.
\end{equation}
Here $\Gamma_{\beta\gamma}^\mu$ denotes the
Christoffel symbols and $\g^{\beta\gamma}$ the
inverse matrix of $ \g_{\beta\gamma}$.
Condition (\ref{eq:9}) is known as the {\it harmonic gauge}. Since
(\ref{eq:1}) is a characteristic (see e.g. \cite{DeTurck}), it is
impossible to solve it in the present form. However, under the
harmonic gauge (\ref{eq:9}), the vacuum Einstein equations
(\ref{eq:9}) is equivalent to the {\it reduced Einstein equations:}
\begin{equation}
\label{wave equation}
\g^{\alpha\beta}\pa_\alpha \pa_\beta
\g_{\gamma\delta}=Q_{\gamma\delta}({\g},\pa{\g}).
\end{equation}
Since $\g_{\alpha\beta}$ has a Lorentzian signature, this is a system
of semi-linear wave equations. The expressions
$Q_{\gamma\delta}({\g},\pa{\g})$ are quadratic functions of the
semi-metric $\g_{\alpha\beta}$ and its first order partial derivatives
$\pa \g_{\alpha\beta}$.

It is well known that the initial data $(\h_{ab},\K_{ab})$ cannot be
prescribed arbitrarily and, in fact, the
Codazzi equations, which relate the curvature of the manifold
$(V,{\g}_{\alpha\beta})$ to the one of $(M,{\h}_{ab})$, lead to the
Einstein constraint equations:
\begin{equation}
\label{eq:3}
\left\{\begin{array}{l}
 R({\bf h})-\|\K_{ab}\|^2+({\rm Tr}_{\h} \K_{ab})^2 =0  \qquad( {\rm
Hamilton \ constraint}) \\ D_a \K_b^{a}- D_b({\rm Tr}_{\h}
\K_{ab})=0  \qquad   \qquad\  \  \   ( {\rm momentum \ constraint})
\end{array}\right.
\end{equation}
Here $R(\h)$ and $D_a$ are the scalar and the covariant derivative
with respect to the metric $\h_{ab}$.

{\it
 The fulfillment of the constraint equations (\ref{eq:3}) is a
necessary condition  for the solution of the wave equations (\ref{wave
equation}) with the initial data (\ref{eq:2}) to satisfy the vacuum
Einstein vacuum equations (\ref{eq:1})}.  We refer
to \cite{bartnik04}, \cite{Rendall_book}, \cite{HAE} or
\cite{wald84:_gener_relat}  for a discussion of  this
fact.

\noindent\textbf{Conclusion:}
We conclude that
any solution of the Cauchy problem for the Einstein
 equations (\ref{eq:1}) includes the treatment of the following  two
problems:
\begin{enumerate}
\item[(i)] Solutions to  constraints (\ref{eq:3}), which can be
reduced to  elliptic equations;
\item[(ii)]  Solutions to the reduced Einstein  equations (\ref{wave
equation}) with  the initial data (\ref{eq:2}).
\end{enumerate}

We will deal with  these two problems in asymptotically flat manifolds
(AF). A Riemanian 3-manifold $(M,{\h_{ab}})$ is  (AF) if there is a
compact subset $K$ such that $M\setminus K$ is diffeomorphic to $
\mathbb{R}^3\setminus B_1(0)$ and the metric $\h_{ab}$ tends to the
identity $\e_{ab}$ at infinity.

\subsection{Existence of the evolution equations}
\label{sec:Existence of the evolution equations}

The unknowns  $\g_{\alpha\beta}$ are functions of $t=x^0$ and
$x^a$, $a=1,2,3$. We assume that for $t=0$, the initial data
(\ref{eq:2}) are given on
a space-like (AF) hypersurface $M\simeq \mathbb{R}^3$ and denote the
Bessel potential space on $\mathbb{R}^3$ by $H^s$.

The following classical result was first proved by  Y. Choquet-Bruhat
\cite{choquet--bruhat52} for $s\geq 3$ and improved by T. Hughes, T.
Kato and J. Marsden \cite{hughes76:_well}.
\vskip 5mm
\noindent
{\bf Theorem A.}
{\it
Assume  the following hold:
 \begin{equation}
\label{eq:6}
({\g}_{\alpha\beta}(0)-{\bf m}_{\alpha\beta})\in
H^{s+1}(\mathbb{R}^3), \quad
\pa_t{\g_{\alpha\beta}}(0)\in H^{s}(\mathbb{R}^3);
 \end{equation}

\begin{equation}
 \label{eq:7}
\sup_{|x|=r}|g_{\alpha\beta}(0,x^1,x^2,x^3)-\mathbf{\bf
m}_{\alpha\beta}|\to
0\quad \text{as } r\to\infty,
\end{equation}
 here $|x|=\sqrt{x_1^2+x_2^2+x_3^2}$ and
${\bf m}_{\alpha\beta}$ is the Minkowski metric.

Then for  $s>\frac 3 2$,  there is a positive  $T>0$ and a unique
semi-metric ${\g_{\alpha\beta}}(t)$ which satisfies (\ref{wave
equation}) and such
that
\begin{equation}
\label{eq:10}
({\g_{\alpha\beta}}(t)-{\bf m}_{\alpha\beta})\in C([0,T], H^{s+1})\ \
\text{and}\ \
\pa_t{\g}(t)\in C([0,T], H^{s}).
\end{equation}
In addition,  if the pair \begin{math} \label{eq:8}
({\g_{\alpha\beta}}(0),\pa_t{\g_{\alpha\beta}}(0))
\end{math} satisfyies the constraint the equations (\ref{eq:3}), then
the metric $\g_{\alpha\beta}(t)$ is a solution to the vacuum Einstein
equation
(\ref{eq:1}).}

S. Klainerman and I. Rodnianski \cite{KR1} succeeded in  improving
regularity below the critical index $\frac 3 2$.

\vskip 5mm
\noindent
{\bf Theorem B.}
{\it Assume the conditions (\ref{eq:6}) and  (\ref{eq:7})
of Theorem A hold and $\g_{\alpha\beta}(t)$ is a classical solution to
(\ref{wave
equation}) such that \begin{math}
({\g_{\alpha\beta}}(0),\pa_t{\g_{\alpha\beta}}(0))\end{math} satisfy
the constraint the equations
(\ref{eq:3}). Then for $s>1$ there is a positive  $T>0$ depending on
$\|\pa \g_{\alpha\beta}(0)\|_{H^{s}}$  such that
\begin{equation}
\label{eq:11}
\|\pa{\g_{\alpha\beta}}(t)\|_{L^2_t L^\infty_x}\leq
C\|\pa{\g_{\alpha\beta}}(0)\|_{H^{s}}.
\end{equation} }

\subsection{Solutions of the constraint equations}
\label{sec:Solutions of the constraint equations}

The $H^s$ spaces are  inappropriate for solutions of the constraint
 equations, roughly speaking, because in these spaces the Laplacian is
not invertible.  It turns out that the
Nirenberg-Walker-Cantor weighted Sobolev spaces
 $H_{m,\delta}$
(\cite{nirenberg73:_null_spaces_ellip_differ_operat_r},\cite{
cantor75:_spaces_funct_condit_r}) are suitable for asymptotically flat
manifolds and indeed these spaces have been widely used in General
Relativity. We denote the norm of the  weighted Sobolev spaces by
\begin{equation}
\label{eq:12}
\|u\|_{H_{m,\delta}}^2=\sum_{|\alpha|\leq
m}\|(1+|x|)^{(\delta+|\alpha|)}\pa^\alpha
u\|_{L^2(\mathbb{R}^3)}^2,\quad -\infty<\delta<\infty,
\end{equation}
and the space $H_{m,\delta}$ is the  completion of
$C_0^\infty(\mathbb{R}^3)$ under the norm (\ref{eq:12}).

We may express now the
asymptotically flat condition  for a Riemannian metric $\h_{ab}$ 
 by means of $(\h_{ab}-\e_{ab})\in H_{m,\delta}$.
In this paper the identity is denoted by $\e_{ab}$.

There is an extensive literature for the solutions of the constraint
equations in asymptotically flat manifolds. The following  theorem
was 
proved under various assumptions in \cite{bartnik04},
 \cite{cantor79},\cite{CHY}\cite{
choquet--bruhat81:_ellip_system_h_spaces_manif_euclid_infin},
 \cite{y.00:_einst_euclid}, \cite{OMC},
\cite{maxwell06:_rough_einst}.

\vskip 5mm
\noindent
{\bf Theorem C.}
{\it
Let $m$ be an integer greater or equal to one, $-\frac 3 2 <\delta
<-\frac 1 2$. Given a set of free data $(\bar{\h}_{ab},\bar{\K}_{ab})$
such that
$((\bar{\h}_{ab}-\e_{ab}),\bar{\K}_{ab})\in H_{m+1,\delta}\times
H_{m,\delta+1} $.
Then there are  exists a conformally equivalent data $({\h}_{ab},
\K_{ab})$
which satisfies the constraint equations (\ref{eq:3}). Moreover, there
is a constant $C$ such that
\begin{equation*}
 \left\|\left( \h_{ab}-\e_{ab},\K_{ab}
\right)\right\|_{H_{m+1,\delta}\times
H_{m,\delta+1}}\leq C \left\|\left(
\bar{\h}_{ab}-\e_{ab},\bar{\K}_{ab}
\right)\right\|_{H_{m+1,\delta}\times
H_{m,\delta+1}}.
\end{equation*} }

\subsection{The main result}
\label{sec:The main result}

The disadvantage of the present situation is the inconsistency of the
  Sobolev spaces of  Theorems A and B with those of the constraint
equations. The initial data for the semi-linear wave equations
(\ref{wave equation})  are given in  $H^s$-spaces, while Theorem C
provides {\color{red}{ the}} initial data (solutions to the constraint
equations
(\ref{eq:3})) in $H_{m,\delta}$.  Therefore it is impossible to obtain
a solution to the Cauchy problem for the vacuum Einstein
equations with
initial data which are given in one type of Sobolev spaces.
{\it Our goal is to unify the Sobolev spaces of the constraint and the
evolution equations.}

Before stating the main theorem we need to introduce the extension of
the spaces $H_{m,\delta}$ into fractional order. We denote a
scaling with $\epsilon$ by $f_\epsilon(x)=f(\epsilon x)$.

\begin{defn}[$H_{s,\delta}$ Sobolev spaces]
\label{def:1}
 For $s\geq 0$ and $-\infty<\delta<\infty$, we define the
 $H_{s,\delta}$ norm by
  \begin{equation}
    \label{eq:const:11}
    \left(\|u\|_{H_{s,\delta}}\right)^2=
    \sum_j 2^{( \frac{3}{2} + \delta)2j} \left\| (\psi_j u)_{2^j}
\right\|_{H^{s}}^{2}.
  \end{equation}
The sequence $\{\psi_j\}\subset C_0^\infty(\mathbb{R}^3)$ satisfies
the following: $ \psi_j(x)=1$ on $K_j=\{x:2^{j-3}\leq |x|\leq
2^{j+2}\}$, $j=1,2,...$, $K_0=\{x:  |x|\leq4\}$; 
${\rm sup}(\psi_j)\subset \{x:2^{j-4}\leq |x|\leq 2^{j+3}\}$,
for $j\geq 1$, ${\rm sup}(\psi_0)\subset\{x: |x|\leq 2^{3}\}$;
\begin{math}
 |\partial^\alpha  \psi_j(x)|\leq  C_\alpha
2^{-|\alpha|j}
\end{math}, where the constant $C_\alpha$ does not depend on $j$.

 The space  $H_{s,\delta}$ is the set of all temperate
  distributions having a finite norm given by (\ref{eq:const:11}).
\end{defn}

Triebel \cite{triebel76:_spaces_kudrj2} proved that whenever $s$ is
equal to an integer
$m$, then
\begin{equation}
\label{eq:18}
\sum_{j=0}^\infty 2^{( \frac{3}{2} + \delta)2j} \| (\psi_j  u)_{2^j}
    \|_{H^{m}}^{2} \sim \sum_{|\alpha|\leq
m}\|(1+|x|)^{(\delta+|\alpha|)}\pa^\alpha u\|_{L^2(\mathbb{R}^3)}^2.
\end{equation}
Thus whenever the parameter $s$ is an integer, the norms (\ref{eq:12})
and (\ref{eq:const:11}) are equivalent.

\begin{thm}[Main results]
\label{thm:main}
Let  $s>\frac 3 2$ and $  -\frac 3 2<\delta<-\frac 1 2$.
Given a set of free data \ $(\bar{\h}_{ab},\bar{\K}_{ab})$ such
that $((\bar{\h}_{ab}-\e_{ab}),\bar{\K}_{ab})\in H_{s+1,\delta}\times
H_{s,\delta+1}
$.

\begin{itemize}
 \item[{\rm (i)}]  Then there  exists a conformally equivalent
data $({\h}_{ab}, \K_{ab})$  which satisfies the constraint equations
(\ref{eq:3}). Moreover $(({\h}_{ab}-\e_{ab}),{\K}_{ab})\in
H_{s+1,\delta}\times
H_{s,\delta+1} $ and depend continuously on the norms of
$((\bar{\h}_{ab}-\e_{ab}),\bar{\K}_{ab})$.

\item[{\rm (ii)}] Then there  exists a $T>0 $ and a
semi-metric ${\g}_{\alpha\beta}(t)$  solution to the vacuum  Einstein
equations
(\ref{eq:1})
such that
\begin{equation}
\label{eq:13}
    ({\g}_{\alpha\beta}(t)-{\bf m_{\alpha\beta}})\in
C([0,T],H_{s+1,\delta})\cap
C^1([0,T],H_{s,\delta+1})
  \end{equation}
and
\begin{equation}
\label{eq:15}
\left. \begin{array}{l}
 \left\|({\g}_{\alpha\beta}(t)-{\bf
m}_{\alpha\beta})\right\|_{H_{s+1,\delta}}\\
\left\|\pa_t{\g}_{\alpha\beta}(t)\right\|_{H_{s,\delta+1}}
\end{array}\right\}\leq
C\left\|\left(
\bar{\h}_{ab}-\e_{ab},\bar{\K}_{ab}
\right)\right\|_{H_{s+1,\delta}\times
H_{s,\delta+1}}.
\end{equation}
for $t\in[0,T]$.
The metric $\g_{\alpha\beta}(t)$ is the unique solution to the reduce
Einstein
(\ref{wave equation}) with initial  data $({\h}_{ab}, \K_{ab})$.
\end{itemize}
 \end{thm}

\begin{rem}[Uniqueness]  Since the Ricci tensor is invariant
under
diffeomorphisms, it is impossible to get a unique solution to the
vacuum Einstein equation. Because if ${\bf R}_{\alpha\beta}(\g)=0$
and $\phi$ is a diffeomorphism, then the pull-back
$\phi^*\g_{\alpha\beta}$ also satisfies (\ref{eq:1}). However, it
can be shown that if two metrics $\g_{\alpha\beta}$ and
$\widetilde{\g}_{\alpha\beta}$ satisfy (\ref{eq:1}) and $(
\g_{\alpha\beta}-{\bf m}_{\alpha\beta}), (
\widetilde{\g}_{\alpha\beta}-{\bf m}_{\alpha\beta})\in
H_{s+1,\delta}$, then there is a coordinates transformation
$x^\alpha\to y^\alpha=f^\alpha(x^\mu)$, which preserve the harmonic
condition (\ref{eq:9}) and  such that
$\widetilde{\g}_{\mu\nu}(y)={\g}_{\alpha\beta}\left(x(y)\right)\frac{
\pa x^\alpha}{\pa y^\mu} \frac{\pa x^\beta}{\pa y^\nu}$. This
transformation can be established by means of solutions to linear
wave equations (see \cite{FMA}, \cite{CHY}) with coefficients in the
$H_{s,\delta}$-spaces.  Thus we can apply the
energy estimate
Lemma \ref{lem:Energy estimates} and the tools of Section
\ref{sec:Weighted Sobelev spaces of fractional order} to establish
the
existence of this transformation in $H_{s+1,\delta}$. Previously this
procedure has been applied with one more degree of differentiability,
but recently, Planchon and Rodnianski
found a trick which
allows the obtaining of  the diffeomorphisms without losing
regularity, see
Section 4 in the monograph \cite{CGP_08} for
details. Thus we conclude that
for  asymptotically flats metrics which preserve the harmonic
condition
(\ref{eq:9}) the uniqueness holds up to a diffeomorphism.

\end{rem}

\begin{rem}
We would like to mention that the results of Christodoulou
\cite{Christodoulou_1981} and  Christodoulou
and O'Murchadha \cite{OMC} differ from ours. They assume 
\begin{math}
 \left(\bar\h_{ab}-\e_{ab},\bar\K_{ab} \right)\\ \in
H_{s+1,\delta+\frac 1 2}\times H_{s+1,\delta+\frac 3 2}
\end{math},
while the solutions $\left(
\g_{\alpha\beta}(t)-{\bf m}_{\alpha\beta}
\right)$ belong to $  H_{s+1,\delta}(\Omega_\theta)$, where
$\Omega_\theta$ is
a certain  unbounded  region of $\mathbb{R}^4$. Thus in their
setting, the rates of  fall-off of the initial data and the
semi-metric are different.  In addition, they require  the regularity
condition
$s\geq 3$.

\end{rem}

The idea to solve both the evolution and the constraint equations in
the  weighted Sobolev spaces of fractional order $H_{s,\delta}$  has
previously  appeared in \cite{ICH12} and  \cite{BK5}, but for the
Einstein-Euler systems. The regularity condition for these systems
is higher since they are coupled with a fluid.

The plan of the paper is as follows. In Section \ref{sec:Weighted
Sobelev spaces of fractional order} we present several properties of
the fractional weighted Sobolev spaces. Section \ref{sec:First order
hyperbolic symmetric systems} deals with the reduction of the wave
equations into a first order symmetric hyperbolic systems. The
specific form of these hyperbolic system has an essential role in our
approach.
The energy estimates are established in Section \ref{sec:Energy
Estimates}. In Section \ref{sec:Local}   we treat the existence,
uniqueness and continuity of  semi-linear  first order symmetric
hyperbolic systems in the $H_{s,\delta}$-spaces and the main result is
proved in Section \ref{sec:main}. In this paper Greek indices will
take the values
$0,1,2,3$ while Latin indices $1,2,3$.

\section{Weighted Sobelev spaces of fractional order}
\label{sec:Weighted Sobelev spaces of fractional order}

Here we present  the basic properties of these spaces and the
equivalence between various  norms. All  these
results  were established in the appendices of \cite{BK4} and
\cite{BK5}. At the end of the section we define a norm on product
spaces.

\begin{defn}[Definitions of norms]
\label{dafn:Definitions of norms}
$\quad$
\begin{itemize}
 \item Let $\{\psi_j\}$ be the sequence of functions in Definition
\ref{def:1}.  For any positive $\gamma$ we set
\begin{equation}
    \label{eq:const:14}
{  \|u\|_{H_{s,\delta,\gamma}}^2 = \sum_j 2^{(
\frac{3}{2} + \delta)2j} \| (\psi_j^\gamma u)_{2^j}
\|_{H^{s}}^{2} }
  \end{equation}
and we will use the convention $ \|u\|_{H_{s,\delta,1}}=
\|u\|_{H_{s,\delta}}$.  The subscripts $2^j$ mean a scaling with
$2^j$,  that
is, $ (\psi_j^\gamma u)_{2^j}(x)=(\psi_j^\gamma u)(2^j x)$.

\item
For a non-negative integer $m$ and $\beta\in \mathbb{R}$, the space
$C_\beta^m$ is the set of all functions having continuous partial
derivatives up to order $m$ and such that the norm (\ref{eq:3.1}) is
finite:
\begin{equation}
\label{eq:3.1}
  \|u\|_{C^m_\beta}=\sum_{|\alpha|\leq
m}\sup_{\mathbb{R}^3}\left((1+|x|)^{\beta+|\alpha|} |\pa^\alpha
u(x)|\right).
\end{equation}

\end{itemize}

\end{defn}

\subsection{Some Properties of $H_{s,\delta}$}
\label{sec:Properties}

\begin{prop}\label{prop:2.1}  \quad\newline
\begin{enumerate}
 \item
\label{equivalence:1} {\bf Equivalence of norms $H_{s,\delta}$ and
$H_{s,\delta,\gamma}$:} For any positive $\gamma$,
\begin{equation}
    \label{eq:const:13}
 \|u\|_{H_{s,\delta}}^2=\sum_j 2^{(
\frac{3}{2} + \delta)2j} \| (\psi_j u)_{2^j}
\|_{H^{s}}^{2} \simeq  \sum_j 2^{(
\frac{3}{2} + \delta)2j} \| (\psi_j^\gamma u)_{2^j}
\|_{H^{s}}^{2}= \|u\|_{H_{s,\delta,\gamma}}^2.
  \end{equation}

\item\label{equivalence} {\bf Equivalence of norms (\ref{eq:const:14})
and
(\ref{eq:12}):}  For any nonnegative integer $m$, positive $\gamma$
and $\delta$ there holds
\begin{equation}
 \label{eq:20}
{  \|u\|_{H_{m,\delta,\gamma}}^2 = \sum_j 2^{(
\frac{3}{2} + \delta)2j} \| (\psi_j^\gamma u)_{2^j}
\|_{H^{m}}^{2} }\sim \sum_{|\alpha|\leq
m}\|(1+|x|)^{(\delta+|\alpha|)}\pa^\alpha
u\|_{L^2(\mathbb{R}^3)}^2.
\end{equation}

\item
\label{embedding 1}
{\bf $H_{s,\delta}$-norm of a derivative:}
\begin{equation}
\label{eq:const:17}
 \|\pa_i u\|_{H_{s-1,\delta+1}}\leq \| u\|_{H_{s,\delta}}.
\end{equation}

\item
\label{Algebra} {\bf Algebra:} If $s_1,s_2\geq s$,
$s_1+s_2>s+\frac{3}{2}$\ and
    $\delta_1+\delta_2\geq\delta-\frac{3}{2}$, then
    \begin{equation}
    \label{eq:elliptic part:12}
      \|uv\|_{H_{s,\delta}}\leq C \|u\|_{H_{s_1,\delta_1}}\
      \|v\|_{H_{s_2,\delta_2}}.
    \end{equation}

\item
  \label{Compact embedding}
{\bf Compact embedding:}
Let $0\leq s'<s$ and $\delta'<\delta$, then the embedding
\begin{equation}
\label{eq:elliptic part:8} H_{s,\delta}\hookrightarrow
H_{s',\delta'}.
\end{equation}
is compact.

\item
\label{Embedding}
{\bf Embedding into the continuous:}  If $s>\frac{3}{2} + m$
and $\delta+\frac 3  2\geq\beta$, then
\begin{equation}
\label{eq:em:2} \| {u}\|_{C^m_\beta}\leq C \|{u}\|_{H_{s,\delta}}.
\end{equation}

\item
\label{Moser}
{\bf Third Moser's inequality:} Let $F:\setR^m\to\setR^l$ be
$C^{N+1}$-function such that
 $F(0)\in H_{s,\delta}$  and where  $N\geq [s]+1$. Then there is a
constant $C$ such that for any $u\in H_{s,\delta}$
\begin{equation}
\label{eq:14}
 \|{{F}}(u)\|_{H_{s,\delta}}\leq
C\|{{F}}\|_{C^{N+1}}\left(1+\|u\|_{L^\infty}^N\right)
\|u\|_{H_{s,\delta}}+
 \|{{F}}(0)\|_{H_{s,\delta}}.
\end{equation}
In particular, if $s>\frac 3 2$, then $ u\in L^\infty$ and
$\| F(u)\|_{H_{s,\delta}}\leq C
\|u\|_{H_{s,\delta}}+\|{{F}}(0)\|_{H_{s,\delta}}$.
\item
\label{Difference estimate}
{\bf Difference estimate:} Suppose $F$ is a $C^{N+2}$-function and
$u,v \in H_{s,\delta}\cap
L^\infty$. Then
\begin{equation}
\label{eq:Moser:6}
 \|{{F}}(u)-{{F}}(v)\|_{H_{s,\delta}}\leq
C(\|u\|_{L^\infty},\|v\|_{L^\infty})
 \left(\|u\|_{H_{s,\delta}}+\|v\|_{H_{s,\delta}}\right)
\left\|u-v\right\|_{H_{s,\delta}}.
\end{equation}

\item
\label{Density}
{\bf Density:}
\begin{enumerate}
\item The class $C_0^\infty(\mathbb{R}^3)$ is dense in
$H_{s,\delta}$.

\item Given $u\in H_{s,\delta}$, $s'>s\geq 0$ and $\delta'\geq\delta$.
Then
for $\rho>0$  there is $u_\rho\in C_0^\infty(\mathbb{R}^3)$ and  a
positive constant $C(\rho)$ such that
\begin{equation}
\|u_\rho -u\|_{H_{s,\delta}}\leq \rho\quad \text{and}\quad
\|u_\rho\|_{H_{s',\delta'}}\leq C(\rho) \|u\|_{H_{s,\delta}}.
\end{equation}
\end{enumerate}

 \item\label{mixed norm}
{\bf Mixed norm estimate:} If $u\in H_{s,\delta} $ and
 $\pa_x u\in H_{s,\delta+1}$, then
\begin{equation}
 \label{eq:3.3}
\|u\|_{H_{s+1,\delta}} \lesssim
\left(\|u\|_{H_{s,\delta}}+\|\pa_x u\|_{H_{s,\delta+1}}\right).
\end{equation}

\end{enumerate}
The proof of (\ref{eq:3.3}) follows from the integral
representation  of the norm (\ref{eq:const:11}) (see \cite{BK4},
\cite{BK5} in the Appendix).
\end{prop}

The density property $(b)$ where proved in \cite{BK5}, \cite{BK4} for
$\delta'=\delta$.  Only a slight modification is needed to include
it also for  $\delta'\geq\delta$ and therefore we leave it to the
reader.

\subsection{Product spaces}
\label{sec:Product spaces}

\begin{defn}[Product spaces]
\label{defn:2} We set
$X_{s,\delta}=H_{s,\delta}\times H_{s,\delta+1}\times H_{s,\delta+1}
$, and the norm of  a vector valued function
$V=(v_1,v_2,v_3)\in
X_{s,\delta}$ is defined by
\begin{equation}
\label{eq:4.75}
 \|V\|_{X_{s,\delta}}^2=\|v_1\|_{H_{s,\delta,2}}^2+\|
v_2\|_{H_{s,\delta+1,2} } ^2+\|v_3\|_{H_{s,\delta+1,2}}^2.
\end{equation}
\end{defn}

We will use the following convention: for a vector valued function
$u:\mathbb{R}\times\mathbb{R}^3\to\mathbb{R}^N$, we set $U=(u,\pa_t
u,\pa_x u)$, where $\pa_x u$ denotes the set of all partial
derivatives in the space variable $x\in\mathbb{R}^3$. Thus
$\|U\|_{X_{s,\delta}}^2={\|u\|_{H_{s,\delta,2}}^2+\|
\pa_t u\|_{H_{s,\delta+1,2} } ^2+\|\pa_x u\|_{H_{s,\delta+1,2}}^2}$.

Essential of our approach is the following observation.
\begin{rem}
 It follows from the Mixed norm estimate (\ref{eq:3.3}) above that
if $U(t,\cdot)\in X_{s,\delta}$,
then $$\|u(t,\cdot)\|_{H_{s+1,\delta}}\lesssim
\|U(t,\cdot)\|_{X_{s,\delta}}.$$
\end{rem}

\section{First order hyperbolic symmetric systems }
\label{sec:First order hyperbolic symmetric systems}

The system of wave equations (\ref{wave equation}) can be transferred
into a first order symmetric hyperbolic system. The specific form of
the hyperbolic system play en important role in our approach.

Letting
\begin{math}
  \h_{\alpha\beta\gamma}=  \partial_\gamma
  \g_{\alpha\beta},
\end{math}
reduces the wave equations (\ref{wave equation})  into
\begin{equation}
  \label{eq:2.3}
  \begin{array}{l}
    \partial_t  \g_{\alpha\beta} = \h_{\alpha\beta  0} \\
    \partial_t  \h_{\gamma\delta  0} = \frac{1}{-\g^{00}}\left\{
    2\g^{0a}\partial_a \h_{\gamma\delta   0} + \g^{ab}\partial_a
\h_{\gamma\delta b} +
C_{\gamma\delta\alpha\beta\rho\sigma}^{
\epsilon\zeta\eta\kappa\lambda\mu}
    \h_{\epsilon\zeta\eta}\h_{\kappa\lambda\mu}
    \g^{\alpha\beta}  \g^{\rho\sigma}\right\}
 \\
   (-\g^{00})^{-1} \g^{ab}\partial_t  \h_{\gamma\delta  a} 
= (-\g^{00})^{-1}\g^{ab}\partial_a
\h_{\gamma\delta 0},
\end{array}
\end{equation}
where the objects
$C_{\gamma\delta\alpha\beta\rho\sigma}^{
\epsilon\zeta\eta\kappa\lambda\mu}$
are a combination of Kronecker deltas with integer coefficients.

We would like now to write system (\ref{eq:2.3}) in a matrix form. We
set
\begin{equation*}
\tilde{g}^{\alpha\beta}=(-\g^{00})^{-1}\g^{\alpha\beta},
\end{equation*}
 where $\g^{\alpha\beta}$ denotes the inverse matrix of
$\g_{\alpha\beta}$.
By introducing the auxiliary vector valued functions
\begin{equation*}
 U=\left(\begin{array}{c}u\\ \pa_t u\\ \pa_x u\end{array}\right)=
\left(\begin{array}{c}\g_{\alpha\beta}-{\bf m}_{\alpha\beta}\\ \pa_t
\g_{\alpha\beta}\\
\pa_x \g_{\alpha\beta}\end{array}\right)=
\left(\begin{array}{c}\g_{\alpha\beta}-{\bf m}_{\alpha\beta}\\
\h_{\alpha\beta0}\\
 \h_{\alpha\beta a}\end{array}\right), \quad a=1,2,3,
\end{equation*}
we can write the system (\ref{eq:2.3}) in the form
\begin{equation}
\label{eq:2.1}
 {\mathcal{A}}^0(u)\pa_t U=\sum_{a=1}^3\left(
{\mathcal{A}}^\alpha(u)+{\mathcal{C}}^a\right) \pa_a
U+{\mathcal{B}}(U)U,
\end{equation}
where
\begin{equation}
 {\mathcal{A}}^0(u)=\left(%
\begin{array}{ccccc}
 \e_{10} & {\bf 0}_{10} & {\bf 0}_{10} & {\bf 0}_{10} & {\bf 0}_{10}
\\
 {\bf 0}_{10} & \e_{10} & {\bf 0}_{10} & {\bf 0}_{10} & {\bf
0}_{10} \\
 {\bf 0}_{10} & {\bf 0}_{10} & \tilde{g}^{11}\e_{10}
&\tilde{g}^{12}\e_{10} &
\tilde{g}^{13}\e_{10}\\
 {\bf 0}_{10} & {\bf 0}_{10} & \tilde{g}^{21}\e_{10} &
\tilde{g}^{22}\e_{10} & \tilde{g}^{23}\e_{10} \\
 {\bf 0}_{10} & {\bf 0}_{10} & \tilde{g}^{31}\e_{10} &
\tilde{g}^{32}\e_{10} &
\tilde{g}^{33}\e_{10} \\
\end{array}%
\right),
\end{equation}

\begin{equation}\label{matrix2}
 {\mathcal{A}}^a(u)=\left(%
\begin{array}{ccccc}
 {\bf 0}_{10} & {\bf 0}_{10} & {\bf 0}_{10} & {\bf 0}_{10} & {\bf
0}_{10} \\
 {\bf 0}_{10} & 2\tilde{g}^{a0}\e_{10}
& (\tilde{g}^{a1}-\delta^{a1})\e_{10} &
(\tilde{g}^{a2}-\delta^{a2})\e_{10} &
(\tilde{g}^{a3}-\delta^{a3})\e_{10} \\
 {\bf 0}_{10} & (\tilde{g}^{a1}-\delta^{a1})\e_{10} & {\bf 0}_{10} &
{\bf
0}_{10} & {\bf 0}_{10}\\
 {\bf 0}_{10} & (\tilde{g}^{a2}-\delta^{a2})\e_{10} & {\bf 0}_{10} &
{\bf
0}_{10} & {\bf 0}_{10} \\
 {\bf 0}_{10} & (\tilde{g}^{a3}-\delta^{a3})\e_{10} & {\bf 0}_{10} &
{\bf
0}_{10}  & {\bf 0}_{10} \\
\end{array}%
\right),
\end{equation}
\begin{equation}\label{matrix4}
 {\mathcal{C}}^a=\left(%
\begin{array}{ccccc}
 {\bf 0}_{10} & {\bf 0}_{10} & {\bf 0}_{10} & {\bf 0}_{10} & {\bf
0}_{10} \\
 {\bf 0}_{10} & {\bf 0}_{10} &
\delta^{a1}\e_{10} &
\delta^{a2}\e_{10} &
\delta^{a3}\e_{10} \\
 {\bf 0}_{10} & \delta^{a1}\e_{10} & {\bf 0}_{10} & {\bf
0}_{10} & {\bf 0}_{10}\\
 {\bf 0}_{10} & \delta^{a2}\e_{10} & {\bf 0}_{10} & {\bf
0}_{10} & {\bf 0}_{10} \\
 {\bf 0}_{10} & \delta^{a3}\e_{10} & {\bf 0}_{10} & {\bf
0}_{10}  & {\bf 0}_{10} \\
\end{array}%
\right),
\end{equation}

and
\begin{equation}\label{matrix3}
 {\mathcal{B}}(U)=\left(%
\begin{array}{ccccc}
 {\bf 0}_{10} & {\e}_{10} & {\bf 0}_{10} & {\bf 0}_{10} & {\bf 0}_{10}
\\
 {\bf 0}_{10} & C^{\kappa\lambda0} & C^{\kappa\lambda1} &
C^{\kappa\lambda2} & C^{\kappa\lambda3} \\
 {\bf 0}_{10} & {\bf 0}_{10} & {\bf 0}_{10} & {\bf 0}_{10} & {\bf
0}_{10}\\
 {\bf 0}_{10} & {\bf 0}_{10} & {\bf 0}_{10} & {\bf 0}_{10} & {\bf
0}_{10} \\
 {\bf 0}_{10} & {\bf 0}_{10} & {\bf 0}_{10} & {\bf 0}_{10}  & {\bf
0}_{10} \\
\end{array}%
\right).
\end{equation}
Here $C^{\kappa\lambda\mu}=(-\g^{00})^{-1}C_{
\gamma\delta\alpha\beta\rho\sigma } ^ {
\varepsilon\zeta\eta\kappa\lambda\mu}\pa_\varepsilon\g_{\zeta\eta}\g^{
\alpha\beta}\g^{\rho\sigma}$.

Apart from the facts that  ${\mathcal{A}}^\alpha(u)$ and
$\mathcal{C}^\alpha$ are symmetric matrices and ${\mathcal{A}}^0(u)$
is positive
definite, they hold three additional  properties: (i) the
matrices
${\mathcal{A}}^\alpha(u)$ doest
not depend on the derivatives of $u$; (ii) the coefficients of
$\pa_t u$ in ${\mathcal{A}}^0(u)$ do not depend on $t$; (iii) it
follows
from Moser type estimate \ref{Moser} and Algebra
\ref{Algebra} of Proposition \ref{prop:2.1} that if
$\g_{\alpha\beta}-{\bf m}_{\alpha\beta}\in H_{s+1,\delta}$, then
$(\tilde{g}^{aa}-1)\in H_{s+1,\delta}$ ($a=1,2,3$) and
$\tilde{g}^{\alpha\beta}\in H_{s+1,\delta}$ whenever
$\alpha\not=\beta$. Thus the matrices $({\mathcal{A}}^0(u)-\e),
{\mathcal{A}}^a(u)\in
H_{s+1,\delta}$ whenever $\g_{\alpha\beta}-{\bf m}_{\alpha\beta}\in
H_{s+1,\delta}$, while ${\mathcal{C}}^a$ are a constant matrices.
These facts are  crucial for the energy estimates.

\section{Energy Estimates}
\label{sec:Energy Estimates}

We consider here  energy estimates for a first order linear hyperbolic
system of the form
\begin{equation}
 \label{eq:4.34}
 {\mathcal{A}}^0\pa_t\left(\begin{array}{c} u\\ \pa_t
u\\ \pa_x u
\end{array}\right) =\sum_{a=1}^3\left(
{\mathcal{A}}^a+{\mathcal{C}}^a\right)
\pa_a\left(
\begin{array}{c} u \\ \pa_t u \\ \pa_x u\end{array}\right)
+{\mathcal{B}}\left(\begin{array}{c} u\\ \pa_t u\\ \pa_x u
\end{array}\right) +{\mathcal{F}}.
\end{equation}
Here $u:\mathbb{R}\times\mathbb{R}^3\to\mathbb{R}^N$,
${\mathcal{A}}^\alpha=\left({\bf a}^\alpha_{i,j}\right)_{i,j=1,2,3}$,
${\mathcal{C}}^a=\left({\bf c}^a_{i,j}\right)_{i,j=1,2,3}$
 and 
 ${\mathcal{B}}=({\bf b}_{ij})_{i,j=1,2,3}$ are  $5N\times 5N$ block
matrices having the sizes of their
blocks  according to the following structure
\begin{equation}\label{eq:4.35}
\left( \begin{array}{c|c|c}
N\times N & N\times N & N \times 3N\\ \hline
N\times N & N\times N & N \times 3N\\ \hline
3N\times N & 3N\times N & 3N \times 3N\\
 \end{array}
\right).
\end{equation}

We assume the following conditions:
\begin{subequations}\label{eq:formulation}
  \begin{eqnarray}
&&
        {\bf a}^0_{ij}=0\ \text{ for} \  i\not=j; \ {\bf
a}^0_{11}={\bf a}^0_{22}=\e;
        \label{item:L1}\\
&& {\bf a}_{33}^0 \quad \text{is symmetric and }
 \frac{1}{c_0}v^Tv\leq
v^T{\bf a}_{33}^0 v\leq c_0v^Tv \quad \forall v\in \mathbb{R}^{3N};
\label{eq:4.2}
   \\
&& {\mathcal{A}}^0(t,\cdot)-\e\in H_{s+1,\delta};
        \label{item:L4}\\ &&
\pa_t{\mathcal{A}}^0(t,\cdot)\in L^\infty; \label{item:L8}
\\ && {\mathcal{A}}^a\ \text{ are symmetric}\ \text{ with}\ \
{\bf a}_{i1}^a={\bf a}_{1j}^a={\bf 0}, \ a=1,2,3;
\label{item:L9} \\ &&
{\mathcal{A}}^a(t,\cdot)\in H_{s+1,\delta}, \quad a=1,2,3;
\label{item:L5}
\\ && {\mathcal{C}}^a\ \text{ are constant and  symmetric}\ \text{
with}\
\
{\bf c}_{i1}^a={\bf c}_{1j}^a={\bf 0}, \ a=1,2,3;
\label{item:L2}
\\ &&
 {\bf b}_{i1}={\bf 0} \ \text{and }\ {\bf b}_{1,j} \ \text{
are constant}, \ \ i,j=1,2,3; \label{item:L10}\\ &&
{\widetilde{\mathcal{B}}}(t,\cdot):=({\bf b}_{ij})_{i,j=2,3}\in
H_{s,\delta+1};\label{item:L6}
\\ &&
{\mathcal{F}}(t,\cdot)\in
H_{s,\delta+1}.\label{item:L7}
  \end{eqnarray}
\end{subequations}

Note  that any system which is originated from a linearization of
 (\ref{eq:2.3}) meets the above requirements.

\subsection{$H_{s,\delta}$ - energy estimates}
\label{subsec:H-energy estimates}

We define an inner-product on $X_{s,\delta}$ in accordance with the
equations (\ref{eq:4.34}). Let
\begin{equation*}\La^s(u)=(1-\Delta)^{\frac{s}{2}}(u)={\mathcal{F}}^{
-1}
\left((1+|\xi|^2)^
{  \frac{s}{2}}{\mathcal{F}}\right)(u),
\end{equation*}
where  $\mathcal{F}$ denote the Fourier transform.

\begin{defn}[Inner-product]
\label{def:inner-product}$\quad$
\begin{itemize}
 \item {\bf Inner-product on $L^2$:} For vector valued
functions $f,g\in L^2$, we set
\begin{equation}
 \label{eq:inner4}
\left\langle  f,g\right\rangle_{L^2}= \int f^Tg dx,
\end{equation}
where $f^T$ denotes the transpose matrix.
 \item {\bf Inner-product on $H_{s,\delta}$:} For $v_1,\phi_1\in
H_{s,\delta}$, we set
\begin{equation}
 \label{eq:inner1}
\left\langle  v_1,\phi_1\right\rangle_{s,\delta}= \sum_{j=0}^\infty
2^{( \frac{3}{2} + \delta)2j}\left\langle \Lambda^{s}\left(\psi_j^2
v_1\right)_{2^j},
\Lambda^{s}\left(\psi_j^2\phi_1\right)_{2^j}\right\rangle_{L^2}.
\end{equation}
Recall that the subscripts $2^j$ mean a scaling (see Definition
\ref{dafn:Definitions of norms}).
 \item {\bf A weighted inner-product on $H_{s,\delta+1}\times
H_{s,\delta+1} $:} For a matrix ${\bf a}_{33}^0$ which satisfies
(\ref{eq:4.2}) and $(v_2,v_3),(\phi_2,\phi_3)\in H_{s,\delta+1}\times
H_{s,\delta+1} $, we set {\small
\begin{equation}
 \label{eq:inner2}
\begin{split}
&\left\langle  \left(\begin{array}{c}v_2\\
v_3\end{array}\right),\left(\begin{array}{c}\phi_2\\
\phi_3\end{array}\right)\right\rangle_{s,\delta+1,{\bf a}_{33}^0}
\\  &=\sum_{j=0}^\infty 2^{( \frac{3}{2} + \delta+1)2j}\left\langle
\Lambda^{s}\left(\psi_j^2 \left(\begin{array}{c}v_2\\
v_3\end{array}\right)\right)_{2^j},\left(\begin{array}{cc}
\e & {\bf 0} \\ {\bf 0} & {\bf a}_{33}^0\\ \end{array}\right)_ { 2^j }
\Lambda^{s}\left(\psi_j^2\left(\begin{array}{c}\phi_2\\
\phi_3\end{array}\right)\right)_{2^j}\right\rangle_{L^2}.
\end{split}
\end{equation}}
\end{itemize}
\begin{itemize}
 \item {\bf Inner-product on $X_{s,\delta}$:}
For a matrix ${\mathcal{A}}^0$ which satisfies
(\ref{item:L1})-(\ref{eq:4.2}) and $V=(v_1,v_2,v_3),
\Phi=(\phi_1,\phi_2,\phi_3)\in X_{s,\delta}$, we set
\begin{equation}
 \label{eq:inner3}
 \la V,\Phi\ra_{X_{s,\delta,{\mathcal{A}}^0}}
 =\left\langle  v_1,\phi_1\right\rangle_{s,\delta}+
\left\langle  \left(\begin{array}{c}v_2\\
v_3\end{array}\right),\left(\begin{array}{c}\phi_2\\
\phi_3\end{array}\right)\right\rangle_{s,\delta+1,{\bf a}_{33}^0}
\end{equation}
\end{itemize}

\end{defn}

We denote by $\|V\|_{X_{s,\delta,{\mathcal{A}}^0}}  $ the norm which
is
associated
with the inner-product (\ref{eq:inner3}).

From condition (\ref{eq:4.2}) we see that

\begin{equation}
 \label{eq:4.6}
\begin{split}
\| V\|_{X_{s,\delta,{\mathcal{A}}^0}}^2 & \leq c_0\left\{
\sum_{j=0}^\infty 2^{(
\frac{3}{2} + \delta)2j} \|\left(\psi_j^2
v_1\right)_{2^j}\|_{H^{s}}^2\right.
\\ & + \left.
\sum_{j=0}^\infty 2^{( \frac{3}{2} + \delta+1)2j}\left[
\|\left(\psi_j^2   v_2\right)_{2^j}\|_{H^{s}}^2+ \|\left(\psi_j^2
v_3\right)_{2^j}\|_{H^{s}}^2\right]\right\}
\\ &  =
c_0\left\{
\|v_1\|_{H_{s,\delta,2}}^2+\|v_2\|_{H_{s,\delta+1,2}}^2+\|v_3\|_{H_{s,
\delta+1,2}}^2\right\}
\\ & =
c_0\| V\|_{X_{s,\delta}}^2.
\end{split}
\end{equation}
Thus  we have shown:

\begin{cor}[Equivalence of norms]
 \label{cor:Equivalence of norms}
The norms which are defined by the inner-product (\ref{eq:inner3}) and
(\ref{eq:4.75}), satisfy
\begin{equation}
 \frac 1 {\sqrt{c_0}}\| V\|_{X_{s,\delta}}\leq \|
V\|_{X_{s,\delta,{\mathcal{A}}^0}}\leq \sqrt{c_0} \|
V\|_{X_{s,\delta}}.
\end{equation}
\end{cor}

For a vector valued function $u(t,x)$, we set $u(t)=u(t,x)$,
$U(t)=(u(t),\pa_t u(t),\pa_x u(t))$ and the energy of $U(t)$
is denoted by
\begin{equation}
 \label{eq:4.7}
E(t)=\frac{1}{2}\la U(t),U(t)\ra_{X_{s,\delta,{\mathcal{A}}^0}}.
\end{equation}
The energy estimate in the product
space $X_{s,\delta}$ is indispensable tool of our method. The next
Lemma establishes it and its proof relies on  tedious computations.
The essential point is that fact   that  ${\mathcal{A}}^\alpha\in
H_{s+1,\delta}$. This enables to use the Kato-Ponce commutator
estimate (Theorem \ref{thm:energy-estimates:1}) with the
pseudodifferential operator $\Lambda^s\pa_x$ rather then
$\Lambda^s$ as in \cite{BK4}.

\begin{lem}[Energy estimates]
\label{lem:Energy estimates}
Let $s>\frac 3 2$, $\delta\geq -\frac 3 2$ and assume  the
coefficients   of (\ref{eq:4.34}) satisfy
conditions (\ref{eq:formulation}). If
$U(t,\cdot)\in C_0^\infty(\mathbb{R}^3)$ is a solution to the linear
system (\ref{eq:4.34}), then
\begin{equation}
 \label{eq:4.8}
\frac{d}{dt}E(t)\leq Cc_0 \left( E(t) +1\right),
\end{equation}
 where the constant $C$ depends
on $\|({\mathcal{A}}^0-\e)\|_{H_{s+1,\delta}}$,
$\|{\mathcal{A}}^a\|_{H_{s+1,\delta}}$,
$\|\widetilde{{\mathcal{B}}}\|_{H_{s,\delta+1}}$,
$\|{\mathcal{F}}\|_{H_{s,\delta+1}}$, $\|\pa_t
{\mathcal{A}}^0\|_{L^\infty}$,
$s$ and
$\delta$.
\end{lem}

An essential tool for deriving these estimates is the Kato \&
Ponce commutator estimate \cite{Taylor91}.
\begin{thm}[Kato and Ponce]
  \label{thm:energy-estimates:1}
  Let $P$ be a pseudodifferential operator in the class $S^s_{1,0}$,
$f\in H^s\cap C^1$, $g\in H^{s}\cap L^{\infty}$ and   $s>0$. Then
\begin{equation}
\label{eq:energy-estimates:9} \| P(fg)-f
P(g)\|_{L^2} \leq C \left\{ \|\nabla f\|_{L^{\infty}}
  \|g\|_{H^{s-1}} + \|f\|_{H^s}\|g\|_{L^{\infty}}
\right\}.
\end{equation}
\end{thm}

\noindent
\textbf{Proof} (of Lemma \ref{lem:Energy estimates}){\bf .}
Taking into account the structure of the inner-product
(\ref{eq:inner3}) we see that
\begin{equation}
 \label{eq:4.10} \begin{split} \frac{d}{dt}E(t) &= \left\langle
u,\pa_t u\right\rangle_{s,\delta} +\left\langle
\left(\begin{array}{c}\pa_t u\\ \pa_x
u\end{array}\right),\pa_t\left(\begin{array}{c}\pa_t u\\ \pa_x
u\end{array}\right)\right\rangle_{s,\delta+1,{\bf a}_{33}^0} \\  + &
\frac 1 2\sum_{j=0}^\infty 2^{( \frac{3}{2} + \delta+1)2j}\left\langle
\Lambda^{s}\left(\psi_j^2 \left(\pa_x
u\right)\right)_{2^j},\pa_t\left({\bf a}_{33}^0\right)_{2^j}
\Lambda^{s}\left(\psi_j^2\left( \pa_x
u\right)\right)_{2^j}\right\rangle_{L^2} \end{split}.
\end{equation}
The infinite sum of the right hand side of (\ref{eq:4.10}) is less
than
\begin{equation}
 \label{eq:4.11}
 \sqrt{3N}\|\pa_t\left({\bf a}_{33}^0\right)\|_{L^\infty}
\sum_{j=0}^\infty 2^{( \frac{3}{2}+
\delta+1)2j}\|\left(\psi_j^2 \pa_x
u\right)_{2^j}\|_{H^{s}}^2=\sqrt{3N}
 \|\pa_t{\bf a}_{33}^0\|_{L^\infty}\|\pa_x u\|_{H_{s,\delta+1,2}}^2
\end{equation}
and
\begin{equation}
 \label{eq:4.5}
\left\langle  u,\pa_t u\right\rangle_{s,\delta}\leq
\|u\|_{H_{s,\delta,2}}\|\pa_t u\|_{H_{s,\delta,2}}\leq \frac 1 2
\left(\|u\|_{H_{s,\delta,2}}^2+\|\pa_t
u\|_{H_{s,\delta+1,2}}^2\right).
\end{equation}

We turn now  to difficult task, this is  the estimation of
$\left\langle  \left(\begin{array}{c}\pa_t u\\ \pa_x
u\end{array}\right),\pa_t\left(\begin{array}{c}\pa_t u\\ \pa_x
u\end{array}\right)\right\rangle_{s,\delta+1,{\bf a}_{33}^0}$.
Setting
\begin{equation}
\label{eq:4.36} E_{\pa_t}(j)=\left\langle
\Lambda^{s}\left(\psi_j^2\pa_tu\right)_{2^j},
\Lambda^{s}\left(\psi_j^2\pa_t(\pa_t u)\right)_{2^j}\right\rangle
_{L^2},
\end{equation}
\begin{equation}
\label{eq:4.37}
 E_{\pa_x}(j)=\left\langle
\Lambda^{s}\left(\psi_j^2\pa_xu\right)_{2^j},\left({\bf
a}^0_{33}\right)_{2^j} \Lambda^{s}\left(\psi_j^2\pa_t(\pa_x
u)\right)_{2^j}\right\rangle_{L^2},
\end{equation}
 and using the specific form of (\ref{eq:inner2})  we  see that
\begin{equation}
 \label{eq:4.14}
\left\langle  \left(\begin{array}{c}\pa_t u\\ \pa_x
u\end{array}\right),\pa_t\left(\begin{array}{c}\pa_t u\\ \pa_x
u\end{array}\right)\right\rangle_{s,\delta+1,{\bf a}_{33}^0}=\sum_{j=0
}^\infty 2^{(
\frac{3}{2}+\delta+1)2j}\left[E_{\pa_t}(j)+E_{\pa_x}(j)\right].
\end{equation}

We first estimate $E_{\pa_x}(j)$. For that purpose we define a
sequence of functions
\begin{equation}
 \label{eq:4.12}
\Psi_k(x)=\left( \sum_{j=0}^\infty\psi_j(x)\right)^{-1}\psi_k(x),
\end{equation}
where $\{\psi_j\}$ is the sequence defined in Definition \ref{def:1}.
This sequence has the following properties: $\Psi_k\in
C_0^\infty(\mathbb{R}^3)$, $|\pa_\alpha \Psi_k(x)|\leq C_\alpha
2^{-k}$, $\sum_{k=0}^\infty \Psi_k(x)=1$ and
\begin{equation}
 \label{eq:4.13}
\Psi_k(x)\psi_j(x)\not= 0 \quad\text{only for}\quad k=j-3,...,j+4.
\end{equation}
We will use the convention that  $\Psi_{j-m}\equiv 0$ whenever
$j-m<0$. Then {
\begin{equation}
 \label{eq:4.15}
\begin{split} 
E_{\pa_x}(j) &=\left\langle  \Lambda^{s}\left[
\left(\psi_j^2\pa_xu\right)_{2^j}\right] ,\left({\bf
a}^0_{33}\right)_{2^j}
\Lambda^{s}\left[\left(\sum_{k=0}^\infty\Psi_k\right)_{2^j}
\left(\psi_j^2\pa_t(\pa_x u)\right)_{2^j}\right]\right\rangle_{L^2}
\\  & = 
\sum_{k=j-3}^{j+4} \left\langle  \Lambda^{s}\left[
\left(\psi_j^2\pa_xu\right)_{2^j}\right] ,\left({\bf
a}^0_{33}\right)_{2^j}
\Lambda^{s}\left[\left(\Psi_k\psi_j^2\pa_t(\pa_x
u)\right)_{2^j}\right]\right\rangle_{L^2}\\ &=:\sum_{k=j-3}^{j+4}
E_{\pa_x}(j,k).
\end{split}
\end{equation}}

Our aim now is to take ${\bf a}_{33}^0$ across $\Lambda^{s}$ in
(\ref{eq:4.37}),  and then we can use the fact that $U$ satisfies
equation (\ref{eq:4.34}). In order to do it we will use the commutator
estimate, Theorem \ref{thm:energy-estimates:1}. However, if  use the
commutator (\ref{eq:energy-estimates:9}) directly
with $f=\left(\Psi_k\right)_{2^j}$ and
$g=\left(\Psi_k\psi_j^2\pa_t\pa_x u\right)_{2^j}$, then we  get
$\|\left(\Psi_k\psi_j^2\pa_t\pa_x
u\right)_{2^j}\|_{L^\infty}\lesssim \|\left(\Psi_k\psi_j^2\pa_t\pa_x
u\right)_{2^j}\|_{H^{s}}$ by Sobolev inequality. That would leads
to the condition
$s-1>\frac{3}{2}$, and therefore we would not obtain the desired
regularity. In order to avoid  this,  we  write
\begin{equation}
 \label{eq:4.16}
\left(\Psi_k\psi_j^2\pa_t\pa_x u\right)_{2^j}=\frac{1}{2^j}
\pa_x\left(\Psi_k\psi_j^2\pa_t u\right)_{2^j}  -\left(
\pa_x\left(\Psi_k\psi_j^2\right)\right)_{2^j}\left(\pa_t
u\right)_{2^j},
\end{equation}
then
\begin{equation}
 \label{eq:4.17}
\begin{split}
&\ \ \ \ \Lambda^{s}\left[ \left(\Psi_k\psi_j^2\pa_t\pa_x
u\right)_{2^j}\right]
\\ & = \frac{1}{2^j}
\Lambda^{s}\left[\pa_x\left(\Psi_k\psi_j^2\pa_t u\right) _{2^j}\right]
 -\Lambda^{s}\left[ \left(
\pa_x\left(\Psi_k\psi_j^2\right)\right)_{2^j}\left(\pa_t
u\right)_{2^j}\right]
\\ & =
\frac{1}{2^j} \left\lbrace \left(\Lambda^{s} \pa_x\right)\left[
\left(\Psi_k\psi_j^2\pa_t u\right)_{2^j}\right]
-(\Psi_k)_{2^j}\left(\Lambda^{s}\pa_x\right)\left[ \left(\psi_j^2\pa_t
u\right)_{2^j}\right] \right\rbrace
\\ & +
\frac{1}{2^j}\left[
(\Psi_k)_{2^j}\left(\Lambda^{s}\pa_x\right)\left(\psi_j^2\pa_t
u\right)_{2^j}\right]
\\ & -
 \Lambda^{s}\left[ \left(
\pa_x\left(\Psi_k\psi_j^2\right)\right)_{2^j}\left(\pa_t
u\right)_{2^j}\right].
\end{split}
\end{equation}
Inserting the last three expressions in each term of $
E_{\pa_x}(j,k)$ in the right hand side of (\ref{eq:4.15}) 
results in
\begin{equation}
 \label{eq:4.18}
\begin{split}
&  E_{\pa_x}(j,k)  = \\ & \frac{1}{2^j}\left\langle
\Lambda^{s}\left[
\left(\psi_j^2\pa_xu\right)_{2^j}\right] ,\left({\bf
a}^0_{33}\right)_{2^j}\left\lbrace \left(\Lambda^{s}
\pa_x\right)\left[ \left(\Psi_k\psi_j^2\pa_t
u\right)_{2^j}\right]\right.\right. \\ &
- \left.\left.(\Psi_k)_{2^j}\left(\Lambda^{s}\pa_x\right)\left[
\left(\psi_j^2\pa_t
u\right)_{2^j}\right]\right\rbrace \right\rangle_{L^2}
\\ & +
 \frac{1}{2^j}\left\langle \Lambda^{s}\left[
\left(\psi_j^2\pa_xu\right)_{2^j}\right] ,\left({\bf
a}^0_{33}\right)_{2^j}(\Psi_k)_{2^j}\left(\Lambda^{s}\pa_x\right)\left
[ \left(\psi_j^2\pa_t u\right)_{2^j}\right] \right\rangle_{L^2}
\\ & -
\left\langle\Lambda^{s}\left[ \left(\psi_j^2\pa_xu\right)_{2^j}\right]
,\left({\bf a}^0_{33}\right)_{2^j}\Lambda^{s}\left[ \left(
\pa_x\left(\Psi_k\psi_j^2\right)\right)_{2^j}\left(\pa_t
u\right)_{2^j}\right]\right\rangle_{L^2}
\\ & =:
E_{\pa_x}(a,j,k)+E_{\pa_x}(b,j,k)+E_{\pa_x}(c,j,k).
\end{split}
\end{equation}

\noindent
\textbf{Estimation of ${ E_{\pa_x}(a,j,k)}$:}  Applying the  Cauchy
Schwarz inequality we get
\begin{eqnarray}
 \label{eq:4.19}
\nonumber
|E_{\pa_x}(a,j,k)| &\leq  & \frac{\sqrt{3N}}{2^j}
\left\|\Lambda^{s}\left(\psi_j^2\pa_xu\right)_{2^j}\right\|_{L^2}
\left\|\left( {\bf a}^0_{33}\right)_{2^j}\right\|_{L^\infty}\\ &
\times &
\left\|\left(\Lambda^{s}\pa_x\right)\left(\Psi_k\psi_j^2\pa_t
u\right)_{2^j}-(\Psi_k)_{2^j}\left(\Lambda^{s}
\pa_x\right)\left(\psi_j^2\pa_t u\right)_{2^j}\right\|_{L^2}.
\end{eqnarray}
The advantage of (\ref{eq:4.16}) is that $\Lambda^{s}\pa_x\in
OPS^{s+1}_{1,0}$, hence the Kato-Ponce
 (\ref{eq:energy-estimates:9}) is applied   with $s+1$ rather
than $s$. Therefore
\begin{equation}
 \label{eq:4.20}
\begin{split}
& \left\|\left(\Lambda^{s}\pa_x\right)\left(\Psi_k\psi_j^2\pa_t
u\right)_{2^j}-(\Psi_k)_{2^j}\left(\Lambda^{s}
\pa_x\right)\left(\psi_j^2\pa_t u\right)_{2^j}\right\|_{L^2}\\ \ \ \
\lesssim  &
\left\|\nabla (\Psi_k)_{2^j}\right\|_{L^\infty}\left\|(\psi_j^2\pa_t
u)_{2^j}\right\|_{H^{s}}+\left\|(\Psi_k)_{2^j}\right\|_{H^{s+1}}
\left\| (\psi_j^2\pa_t u)_{2^j}\right\|_{L^\infty}
\end{split}.
\end{equation}
From the properties of $\psi_j$ (see Definition \ref{def:1} ) and
(\ref{eq:4.12}), we see that  $\|\pa^\alpha
(\Psi_k)_{2^j}\|_{L^\infty}\leq C$, where the constant $C$ is
independent of $j$ and $k $. Hence both $\left\|\nabla
(\Psi_k)_{2^j}\right\|_{L^\infty}$ and $\left\|
(\Psi_k)_{2^j}\right\|_{H^s}$ are bounded by a certain constant
independent of $j$ and $k $.
For $s>\frac{3}{2}$ the Sobolev inequality yields $\left\|
(\psi_j^2\pa_t u)_{2^j}\right\|_{L^\infty}\lesssim \left\|
(\psi_j^2\pa_t u)_{2^j}\right\|_{H^{s}}$,  and combining these with
inequality (\ref{eq:4.19}) we get
\begin{equation}
 \label{eq:4.23}
|E_{\pa_x}(a,j,k)|\lesssim \left\|\left( {\bf
a}^0_{33}\right)_{2^j}\right\|_{L^\infty}
\left\|\left(\psi_j^2\pa_xu\right)_{2^j}\right\|_{H^{s}}
\left\|\left(\psi_j^2\pa_tu\right)_{2^j}\right\|_{H^{s}}.
\end{equation}

\noindent
\textbf{Estimation of $E_{\pa_x}(c,j,k)$:} Since
\begin{math}\left(\pa_x(\Psi_k\psi_j^2)\right)_{2^j}(\pa_t
u)_{2^j}=(F_{j,k})_{2^j}(\psi_j\pa_t u)_{2^j} \end{math}  and the
partial derivatives of $ (F_{j,k})_{2^j}$ up to order $[s]$ are
bounded by a constant $C$ independent of $j$ and $k$, there holds
\begin{equation}
 \label{eq:4.24}
\begin{split}
|E_{\pa_x}(c,j,k)| & \lesssim
  \left\|\left( {\bf a}^0_{33}\right)_{2^j}\right\|_{L^\infty}
\left\|\left(\psi_j^2\pa_xu\right)_{2^j}\right\|_{H^{s}}
\left\|(\psi_j\pa_t u)_{2^j}\right\|_{H^{s}}.
\end{split}
\end{equation}

\noindent
\textbf{Estimation of $E_{\pa_x}(b,j,k)$:} We see from
(\ref{eq:4.18}) that in order to use equation (\ref{eq:4.34}) we need
to commute $(\Psi_k{\bf a}_{33}^0)_{2^j}$ with
$\left(\Lambda^{s}\pa_x\right)$.  Therefore we write
\begin{equation}
 \label{eq:4.54}
\begin{split}
& \left(\Psi_k{\bf
a}_{33}^0\right)_{2^j}\left(\Lambda^{s}\pa_x\right)\left[
(\psi_j^2\pa_t u)_{2^j}\right]
\\  = &
\left(\Psi_k{\bf
a}_{33}^0\right)_{2^j}\left(\Lambda^{s}\pa_x\right)\left[
(\psi_j^2\pa_t u)_{2^j}\right] -
\left(\Lambda^{s}\pa_x\right)\left[\left(\Psi_k{\bf
a}_{33}^0\right)_{2^j}(\psi_j^2\pa_t u)_{2^j}\right]
\\ + &
\left(\Lambda^{s}\pa_x\right)\left[\left(\Psi_k{\bf
a}_{33}^0\right)_{2^j}(\psi_j^2\pa_t u)_{2^j}\right]
\\  = &
\left(\Psi_k{\bf
a}_{33}^0\right)_{2^j}\left(\Lambda^{s}\pa_x\right)\left[
(\psi_j^2\pa_t u)_{2^j}\right] -
\left(\Lambda^{s}\pa_x\right)\left[\left(\Psi_k{\bf
a}_{33}^0\right)_{2^j}(\psi_j^2\pa_t u)_{2^j}\right]
\\ + &
\Lambda^{s}\left[\pa_x\left( \Psi_k{\bf a}_{33}^0\psi_j^2\right)_{2^j}
 (\pa_t u)_{2^j}\right]
\\ + &
2^j\Lambda^{s}\left[ \left(\Psi_k\psi_j^2{\bf
a}_{33}^0\right)_{2^j}\left(\pa_t\pa_x u\right)_{2^j}\right].
\end{split}
\end{equation}
Thus $E_{\pa_x}(b,j,k)$ is  sum of three terms:
 $E_{\pa_x}(b,j,k)=E_{\pa_x}(d,j,k)+E_{\pa_x}(e,j,k)+E_{\pa_x}(f,j,
k)$. The first one will be estimated by means Theorem
\ref{thm:energy-estimates:1}, in the second one we use algebra
property of $H^s$ and in the last one brings us to equation
(\ref{eq:4.34}).

We recall that $(\Lambda^{s}\pa_x)\in OPS^{s+1}_{1,0}$, therefore by
Kato-Ponce commutator estimate (\ref{eq:energy-estimates:9}),
\begin{equation}
 \label{eq:4.25}
 \begin{split}
 & \left\|
\left(\Psi_k{\bf
a}_{33}^0\right)_{2^j}\left(\Lambda^{s}\pa_x\right)\left[
(\psi_j^2\pa_t u)_{2^j}\right] -
\left(\Lambda^{s}\pa_x\right)\left[\left(\Psi_k{\bf
a}_{33}^0\right)_{2^j}(\psi_j^2\pa_t u)_{2^j}\right]\right\|_{L^2}
\\  \lesssim &
\left\{\left\|\nabla(\Psi_k {\bf
a}_{33}^0)_{2^j}\right\|_{L^\infty}\left\|(\psi_j^2\pa_t
u)_{2^j}\right\|_{H^{s}}+\left\|(\Psi_k {\bf
a}_{33}^0)_{2^j}\right\|_{H^{s+1}}\left\|(\psi_j^2\pa_t
u)_{2^j}\right\|_{L^\infty}\right\}.
\end{split}
 \end{equation}
From (\ref{eq:4.12}) and (\ref{eq:4.13}) we see that
\begin{equation}
 \label{eq:4.26}
\left\|\nabla(\Psi_k {\bf a}_{33}^0)_{2^j}\right\|_{L^\infty}\leq C_1
\left\|({\bf a}_{33}^0)_{2^j}\right\|_{L^\infty}+C_2
2^j\left\|\nabla{(\bf a}_{33}^0)_{2^j}\right\|_{L^\infty}
\end{equation}
and
\begin{equation}
 \label{eq:4.27}
\begin{split}
& \ \ \ \ \ \ \ \left\|(\Psi_k {\bf
a}_{33}^0)_{2^j}\right\|_{H^{s+1}}
\\ & \leq
C_3 \left\|(\psi_k{\bf a}_{33}^0)_{2^j}\right\|_{H^{s+1}}=C_3
\left\|\left( \psi_k\left( {\bf a}_{33}^0-\e\right)
\right)_{2^j}+\left( \psi_k\e\right)_{2^j}\right\|_{H^{s+1}}
\\  & =
C_3 \left\|\left( \left( \psi_k\left( {\bf a}_{33}^0-\e\right)
\right)_{2^k}\right)_{2^{j-k}}+\left( \left(
\psi_k\e\right)_{k}\right)_{2^{j-k}}\right\|_{H^{s+1}}
\\ & \leq
C_3\left\|\left( \left( \psi_k\left( {\bf a}_{33}^0-\e\right)
\right)_{2^k}\right)_{2^{j-k}}\right\|_{H^{s+1}}+C_3\left\|\left(
\left( \psi_k\e\right)_{k}\right)_{2^{j-k}}\right\|_{H^{s+1}}
\\ & \simeq
C_3\left\lbrace \left\|\left(  \psi_k{\bf
a}_{33}^0-\e\right)_{2^k}\right\|_{H^{s+1}}+1\right\rbrace.
\end{split}
\end{equation}
Thus, the combination of (\ref{eq:4.26}) and (\ref{eq:4.27}) with the
inequality
 $\left\|(\psi_j^2\pa_t u)_{2^j}\right\|_{L^\infty}\lesssim
\left\|(\psi_j^2\pa_t u)_{2^j}\right\|_{H^{s}}$  (keeping in mind
the factor $2^{-j}$ in (\ref{eq:4.18})), leads to
\begin{equation}
\label{eq:4.28}
\begin{split}
&|E_{\pa_x}(d,j,k)| \\ = &\left\langle \Lambda^s \left[ (\psi_j^2\pa_x
u)_{2^j}\right] , \left(\Psi_k{\bf
a}_{33}^0\right)_{2^j}\left(\Lambda^{s}\pa_x\right)\left[
(\psi_j^2\pa_t u)_{2^j}\right]\right. \\  - & \left.
\left(\Lambda^{s}\pa_x\right)\left[\left(\Psi_k{\bf
a}_{33}^0\right)_{2^j}(\psi_j^2\pa_t u)_{2^j}\right]
\right\rangle_{L^2} \\
\lesssim & \left\{\left\|\left( {\bf
a}_{33}^0\right)_{2^j} \right\|_{L^\infty}+\left\|\left( \nabla{\bf
a}_{33}^0\right)_{2^j} \right\|_{L^\infty}+
\left\|\left(  \psi_k({\bf
a}_{33}^0-\e)\right)_{2^k}\right\|_{H^{s+1}}
\right\}
\\   \times &
\left\|(\psi_j^2\pa_x u)_{2^j}\right\|_{H^{s}}\left\|(\psi_j^2\pa_t
u)_{2^j}\right\|_{H^{s}}.
\end{split}
\end{equation}
We turn now to $E_{\pa_x}(e,j,k)=\left\langle \Lambda^s \left[
(\psi_j^2\pa_x
u)_{2^j}\right],\left(\Lambda^{s}\pa_x\right)\left[\left(\Psi_k{\bf
a}_{33}^0\right)_{2^j}(\psi_j^2\pa_t u)_{2^j}\right]
\right\rangle_{L^2}$.  Noting that
\begin{equation}
\begin{split}
& \pa_x\left( \Psi_k{\bf a}_{33}^0\psi_j^2\right)_{2^j}(\pa_t
u)_{2^j}\\ = & \pa_x\left( \Psi_k{\bf
a}_{33}^0\psi_j\right)_{2^j}(\psi_j\pa_t u)_{2^j}+ 2^j\left(
\Psi_k{\bf a}_{33}^0\pa_x\psi_j\right)_{2^j}(\psi_j\pa_t u)_{2^j},
\end{split}
\end{equation}
and applying  the algebra  property of $H^s$, we get
\begin{equation}
 \label{eq:4.29}
\begin{split}
 & \left\|\pa_x\left( \Psi_k{\bf a}_{33}^0\psi_j^2\right)_{2^j}(\pa_t
u)_{2^j}\right\|_{H^{s}}
  \lesssim 
\left\|\pa_x\left( \Psi_k{\bf
a}_{33}^0\psi_j\right)_{2^j}\right\|_{H^{s}}\left\|(\psi_j\pa_t
u)_{2^j}\right\|_{H^{s}}
 \\ &+
2^j\left\|\left( \Psi_k{\bf
a}_{33}^0\pa_x\psi_j\right)_{2^j}\right\|_{H^{s}}\left\|(\psi_j\pa_t
u)_{2^j}\right\|_{H^{s}}.
\end{split}
\end{equation}
Now
\begin{math}
 \label{eq:4.30}
 \left\|\pa_x\left( \Psi_k{\bf
a}_{33}^0\psi_j\right)_{2^j}\right\|_{H^{s}}
\lesssim  \left\|\left( \Psi_k{\bf
a}_{33}^0\right)_{2^j}\right\|_{H^{s+1}}
\end{math}
and
\begin{math}
 \left\|\left( \Psi_k{\bf
a}_{33}^0\pa_x\psi_j\right)_{2^j}\right\|_{H^{s}}
\lesssim  \left\|\left( \Psi_k{\bf
a}_{33}^0\right)_{2^j}\right\|_{H^{s+1}}
\end{math},
hence  by   (\ref{eq:4.29}) and inequality (\ref{eq:4.27}) we get
\begin{equation}
 \label{eq:4.32}
\begin{split}
&|E_{\pa_x}(e,j,k)|\\  \lesssim & \left\{ \left\|\left(\psi_k({\bf
a}_{33}^0-\e)\right)_{2^k}\right\|_{H^{s+1}}+1 \right\}
\left\|(\psi_j^2\pa_x u)_{2^j}\right\|_{H^{s}}\left\|(\psi_j\pa_t
u)_{2^j}\right\|_{H^{s}}.
\end{split}
\end{equation}
In order to use equation (\ref{eq:4.34}) we write
\begin{equation}\label{eq:4.39}
{\mathcal{A}}^a=\left(
\begin{array}{c|cccc}
{\bf 0} &  &   {\bf 0} & \\
\hline
& & &   &\\
{\bf 0} &  &
\widetilde{\mathcal{A}}^a & \\
& &  && \\
\end{array}
\right),\qquad {\mathcal{C}}^a=\left(
\begin{array}{c|cccc}
{\bf 0} &  &   {\bf 0} & \\
\hline
& & &   &\\
{\bf 0} &  &
\widetilde{\mathcal{C}}^a & \\
& &  && \\
\end{array}
\right),\ \ a=1,2,3,
\end{equation}
where $\widetilde{\mathcal{A}}^a=({\bf a}_{ij}^a)_{i,j=2,3}$,
$\widetilde{\mathcal{C}}^a=({\bf c}_{ij}^a)_{i,j=2,3}$ are  symmetric
block matrix and ${\bf c}_{ij}^a$ are constant. Further, let 
 $\{\Psi_k\} $ be the sequence which is defined by (\ref{eq:4.12}),
then
\begin{equation}
 \label{eq:4.38}
\begin{split}
 &\ \ \ \  E_{\pa_t}(j)=\left\langle \Lambda^{s}(\psi_j^2 \pa_t
u)_{2^j},\Lambda^{s}(\psi_j^2 \pa_t(\pa_t u))_{2^j}\right\rangle_{L^2}
\\ & =\left\langle \Lambda^{s}(\psi_j^2 \pa_t
u)_{2^j},\Lambda^{s}\left[\left(\sum_{k=0}^\infty\Psi_k\right)_{2^j}
(\psi_j^2 \pa_t (\pa_t u))_{2^j}\right]\right\rangle_{L^2}
\\  &=
\sum_{k=j-3}^{j+4}\left\langle \Lambda^{s}(\psi_j^2 \pa_t
u)_{2^j},\Lambda^{s}\left[\left(\Psi_k\psi_j^2 \pa_t(\pa_t
u)\right)_{2^j}\right]\right\rangle_{L^2}\\ &=:\sum_{k=j-3}^{j+4}E_{
\pa_t}(j , k).
\end{split}
\end{equation}

From (\ref{eq:4.15}), (\ref{eq:4.18})  and (\ref{eq:4.54})  we
see  that
\begin{equation*}
 E_{\pa_x}(f,j,k)=\left\langle \Lambda^s \left[
(\psi_j^2\pa_x
u)_{2^j}\right],\Lambda^{s}\left[\left(\Psi_k{\bf
a}_{33}^0\right)_{2^j}(\psi_j^2\pa_t \pa_x u)_{2^j}\right]
\right\rangle_{L^2}
\end{equation*}
and since $U(t)$ satisfies (\ref{eq:4.34}), we have obtained
{
\begin{equation}
 \label{eq:4.40}
\begin{split}
& \ \ \  \left\{
E_{\pa_t}(j,k) +
 E_{\pa_x}(f,j,k)
\right\} = \\ &
\left\langle \Lambda^{s}\left[  \left( \psi_j^2\left(\begin{matrix}
 \pa_t u \\
  \pa_x u
\end{matrix}\right)\right)_{2^j}\right] , \Lambda^{s}\left[ \left(
\Psi_k\psi_j^2 \left(\begin{array}{cc}{\e} & {\bf 0}\\ {\bf 0} & {\bf
a}_{33}^0\end{array}\right) \pa_t\left( \begin{matrix}
 \pa_t u \\
  \pa_x u
\end{matrix}\right)\right)_{2^j}\right] \right\rangle_{L^2}=
\\  &  \sum_{a=1}^3
\left\langle \Lambda^{s} \left[ \left( \psi_j^2\left(\begin{matrix}
 \pa_t u \\
  \pa_x u
\end{matrix}\right)\right)_{2^j}\right] , \Lambda^{s}
\left[ \left( \Psi_k\psi_j^2 \left(\begin{matrix}
\left( \widetilde{\mathcal{A}}^a+\widetilde{\mathcal{C}}^a\right)
\pa_a\left(\begin{matrix}
 \pa_t u \\
 \pa_x u
\end{matrix}\right)
\end{matrix}\right)\right)_{2^j}\right] \right\rangle_{L^2}
\\  & + 
\left\langle \Lambda^{s} \left[ \left( \psi_j^2\left(\begin{matrix}
 \pa_t u \\
  \pa_x u
\end{matrix}\right)\right)_{2^j}\right] , \Lambda^{s}\left[ \left(
\Psi_k\psi_j^2 { \left(\begin{matrix}{\bf b}_{22} & {\bf b}_{23}\\
{\bf b}_{32} & {\bf b}_{33}
\end{matrix}\right) \left(\begin{matrix}
 \pa_t u \\
  \pa_x u
\end{matrix}\right)}\right)_{2^j}\right] \right\rangle_{L^2}
\\  & + 
\left\langle \Lambda^{s} \left[ \left( \psi_j^2\left(\begin{matrix}
 \pa_t u \\
  \pa_x u
\end{matrix}\right)\right)_{2^j}\right] , \Lambda^{s}\left[ \left(
\Psi_k\psi_j^2
\left(\begin{matrix} f_2 \\
f_3\end{matrix}\right)\right)_{2^j}\right] \right\rangle_{L^2}.
\end{split}
\end{equation}}
The main difficulty is the estimation of the first term of the right
hand side of (\ref{eq:4.40}). We recall that
$\widetilde{\mathcal{C}}^a$
are
constant and  $\widetilde{\mathcal{A}}^a\in H_{s+1,\delta}$,
therefore we may write
\begin{equation}
 \label{eq:4.41}
\begin{split}
&\left( \Psi_k\psi_j^2\left(
\widetilde{\mathcal{A}}^a+\widetilde{\mathcal{C}}^a\right)
\pa_a \left(\begin{matrix}
 \pa_t u \\
 \pa_x u
\end{matrix}\right)\right)_{2^j} \\ = &\frac{1}{2^j}
\pa_a\left( \Psi_k\psi_j^2\left(
\widetilde{\mathcal{A}}^a+\widetilde{\mathcal{C}}^a\right)
\left(\begin{matrix}
 \pa_t u \\
 \pa_x u
\end{matrix}\right)\right)_{2^j} -
\left( \pa_a\left( \Psi_k\psi_j^2\widetilde{\mathcal{A}}^a\right)
\left(\begin{matrix}
 \pa_t u \\
 \pa_x u
\end{matrix}\right)\right)_{2^j}
\end{split}
\end{equation}
and hence
\begin{equation}
 \label{eq:4.42}
\begin{split}
 & \Lambda^{s}\left[ \left( \Psi_k\psi_j^2\left(
\widetilde{\mathcal{A}}^a+\widetilde{\mathcal{C}}^a\right) \pa_a
\left(\begin{matrix}
 \pa_t u \\
 \pa_x u
\end{matrix}\right)\right)_{2^j}\right]
\\ = & \frac{1}{2^j} \left\{
\left( \Lambda^{s}\pa_a\right) \left[ \left( \Psi_k\left(
\widetilde{\mathcal{A}}^a+\widetilde{\mathcal{C}}^a \right)
\psi_j^2\left(\begin{matrix}
 \pa_t u \\
 \pa_x u
\end{matrix}\right)\right)_{2^j} \right] \right.\\ -& \left.
\left(\Psi_k\left(
\widetilde{\mathcal{A}}^a+\widetilde{\mathcal{C}}^a\right)
\right)_{2^j}\left( \Lambda^{s}
\pa_a\right) \left[ \left( \psi_j^2 \left(\begin{matrix}
 \pa_t u \\
 \pa_x u
\end{matrix}\right)\right)_{2^j}\right]\right\}
\\ +  & \frac{1}{2^j}
 \left(\Psi_k\left(
\widetilde{\mathcal{A}}^a+\widetilde{\mathcal{C}}^a\right)
\right)_{2^j} \Lambda^{s}
 \left[ \pa_a\left( \psi_j^2 \left(\begin{matrix}
 \pa_t u \\
 \pa_x u
\end{matrix}\right)\right)_{2^j}\right]
\\ - &
\frac{1}{2^j}\Lambda^{s}\left[  \pa_a \left(\left(
\Psi_k\psi_j^2\widetilde{\mathcal{A}}^a\right)  \left(\begin{matrix}
 \pa_t u \\
 \pa_x u
\end{matrix}\right)\right)_{2^j}\right].
\end{split}
\end{equation}
The first term of the right hand side of (\ref{eq:4.42}) will be
estimate by Theorem \ref{thm:energy-estimates:1} with
$P=\Lambda^{s}\pa_a$, in the second one the symmetry of
$\widetilde{\mathcal{A}}^a$ will be exploited and  in the third one we
will use
algebra property of $H^s$. In both the first and third we take the
advantage that $\widetilde{\mathcal{A}}^a\in H_{s+1,\delta}$.

\begin{equation*}
 \label{eq:43}
\begin{split}
& \left\|\left( \Lambda^{s}\pa_a\right) \left[
(\Psi_k\left(\widetilde{\mathcal{A}}^a+\widetilde{\mathcal{C}}
^a\right)
)\psi_j^2\left(\begin{matrix}
 \pa_t u \\
 \pa_x u
\end{matrix}\right)\right] -
\left(\Psi_k\left(\widetilde{\mathcal{A}}^a+\widetilde{\mathcal{C}}
^a\right)\right)_ {2^j}
\left( \Lambda^{s}
\pa_a\right) \left[ \left( \psi_j^2 \left(\begin{matrix}
 \pa_t u \\
 \pa_x u
\end{matrix}\right)\right)_{2^j}\right] \right\|_{L^2}
\\ \lesssim &
\left\|\nabla(\Psi_k\left(\widetilde{\mathcal{A}}^a+\widetilde{
\mathcal{C} }
^a\right)
)_{2^j}\right\|_{L^\infty}\left\{\left\| (\psi_j^2\pa_t u
)_{2^j}\right\|_{H^{s}}+\left\| (\psi_j^2\pa_x u
)_{2^j}\right\|_{H^{s}}\right\}
\\   +  &
\left\|(\Psi_k\left(\widetilde{\mathcal{A}}^a+\widetilde{\mathcal{C}}
^a\right)
)_{2^j}\right\|_{H^{s+1}}\left\{\left\| (\psi_j^2\pa_t u
)_{2^j}\right\|_{L^\infty}+\left\| (\psi_j^2\pa_x u
)_{2^j}\right\|_{L^\infty}\right\}.
\end{split}
\end{equation*}

Now, for $k=j-3,...,j+4$,
\begin{equation}
 \label{eq:4.44}
\left\|\nabla\left(\Psi_k\left(\widetilde{\mathcal{A}}^a+\widetilde{
\mathcal{C}} ^a\right)
\right)_{2^j}\right\|_{L^\infty}\lesssim
\left(\left\|\widetilde{\mathcal{A}}^a
\right\|_{L^\infty} +1+ 2^j\left\|\nabla\widetilde{\mathcal{A}}^a
\right\|_{L^\infty}\right),
\end{equation}
and 
\begin{equation}
 \label{eq:4.45}
\begin{split}
&\left\|\left(\Psi_k\left(\widetilde{\mathcal{A}}^a+\widetilde{
\mathcal{C} }
^a\right)\right)_{2^j} \right\|_{H^{s+1}} \lesssim
\left\{\left\|(\psi_k\widetilde{\mathcal{A}}^a
)_{2^j}\right\|_{H^{s+1}}+1\right\} \\ = &\left\{\left\|\left(
(\psi_k\widetilde{\mathcal{A}}^a )_{2^k}\right)_{2^{j-k}}
\right\|_{H^{s+1}} + 1\right\}\lesssim
\left\{\left\|(\psi_k\widetilde{\mathcal{A}}^a
)_{2^k}\right\|_{H^{s+1}}+1\right\}.
\end{split}
\end{equation}
In addition, since $s>\frac 3 2$,
  $$\left\| (\psi_j^2\pa_t u )_{2^j}\right\|_{L^\infty}\lesssim
\left\| (\psi_j^2\pa_t u )_{2^j}\right\|_{H^{s}} \quad \text{
and}\quad \left\| (\psi_j^2\pa_x u )_{2^j}\right\|_{L^\infty}\lesssim
\left\| (\psi_j^2\pa_x u )_{2^j}\right\|_{H^{s}}.$$

Thus,
\begin{equation}
 \label{eq:46}
\begin{split}
\frac{1}{2^j} &\mid \left\langle \Lambda^{s} \left[ \left(
\psi_j^2\left(\begin{matrix}
 \pa_t u \\
  \pa_x u
\end{matrix}\right)\right)_{2^j}\right],
\left( \Lambda^{s}\pa_a\right) \left[
\left(\Psi_k\left(\widetilde{\mathcal{A}}^a+\widetilde{\mathcal{C}}
^a\right) \psi_j^2\left(\begin{matrix}
 \pa_t u \\
 \pa_x u
\end{matrix}\right)\right)_{2^j}\right] \right.
\\ -&
\left.\left(\Psi_k\left(\widetilde{\mathcal{A}}^a+\widetilde{\mathcal{
C}}
^a\right)\right)_{2^j}\left( \Lambda^{s}
\pa_a\right) \left[ \left( \psi_j^2 \left(\begin{matrix}
 \pa_t u \\
 \pa_x u
\end{matrix}\right)\right)_{2^j}\right]\right\rangle_{L^2}\mid
\\ \lesssim  &\left\lbrace \left\|\nabla\widetilde{\mathcal{A}}^a
\right\|_{L^\infty}+\left\|\widetilde{\mathcal{A}}^a
\right\|_{L^\infty} +\left\|(\psi_{k}\widetilde{\mathcal{A}}^a)_{2^k}
\right\|_{H^{s+1}} +1 \right\rbrace \\ \times &
\left\lbrace  \left\| (\psi_j^2\pa_t u )_{2^j}\right\|_{H^{s}}^2
+\left\| (\psi_j^2\pa_x u )_{2^j}\right\|_{H^{s}}^2 \right\rbrace.
\end{split}
\end{equation}
As to the third term of (\ref{eq:4.42}), writing
\begin{equation*}
\pa_a\left(
\Psi_k\psi_j^2\widetilde{\mathcal{A}}^a\right)\left(\begin{matrix}
 \pa_t u \\
 \pa_x u
\end{matrix}\right)=\left( \pa_a(\Psi_k\psi_j\
\widetilde{\mathcal{A}}^a)+2(\Psi_k\widetilde{\mathcal{A}}
^a\pa_a\psi_j\right)\psi_j\left(\begin{matrix}
 \pa_t u \\
 \pa_x u
\end{matrix}\right),
\end{equation*}
and noting  that
\begin{equation*}
\left\|\pa_a\left(\left(
\Psi_k\psi_j\widetilde{\mathcal{A}}^a\right)\right)_{2^j}\right\|_{H^s
}
\leq
\left\|\left(\Psi_k\psi_j\widetilde{\mathcal{A}}^a\right)_{2^j}
\right\|_{
H^{s+1}
}\lesssim
\left\|\left(\psi_k\widetilde{\mathcal{A}}^a\right)_{2^k}\right\|_{H^{
s+1}},
\end{equation*}
then by the embedding
$H^{s+1}\hookrightarrow H^s$ we get that
\begin{equation}
 \label{eq:4.48}
\begin{split}
& \mid \left\langle \Lambda^{s} \left[ \left(
\psi_j^2\left(\begin{matrix}
 \pa_t u \\
  \pa_x u
\end{matrix}\right)\right)_{2^j}\right] ,
\Lambda^{s}\left[\left(  \pa_a\left(
\Psi_k\psi_j^2\widetilde{\mathcal{A}}^a\right)  \left(\begin{matrix}
 \pa_t u \\
 \pa_x u
\end{matrix}\right)\right)_{2^j}\right]  \right\rangle_{L^2}\mid
\\  \lesssim  & \left\lbrace
\left\|(\psi_k\widetilde{\mathcal{A}}^a)_{2^k}\right\|_{H^{s+1}} +
\left\|(\psi_j\widetilde{\mathcal{A}}^a)_{2^j}\right\|_{H^{s+1}}
\right\rbrace \\
\times   &
 \left(\left\| (\psi^2_j\pa_t u)_{2^j}\right\|_{H^{s}}+ \left\|
(\psi_j^2\pa_x u)_{2^j}\right\|_{H^{s}}\right)
\left( \left\| (\psi_j\pa_t u)_{2^j}\right\|_{H^{s}}+ \left\|
(\psi_j\pa_x u)_{2^j}\right\|_{H^{s}}\right).
\end{split}
\end{equation}

We turn now  the second term of (\ref{eq:4.42}).
 Recall $U$ is a $C_0^\infty(\mathbb{R}^3)$, therefore
$\Lambda^{s}(\pa_tu), \Lambda^{s}(\pa_xu)$ are rapidly decreasing
functions. This allows us to make the following operations:
{\small
\begin{equation}
 \label{eq:4.49}
\begin{split}
&\int\pa_a\left\{\left( \Lambda^{s}\left[ \left( \psi_j^2
\left(\begin{matrix}
 \pa_t u \\
  \pa_x u
\end{matrix}\right)\right)_{2^j}\right]\right)^T
\left(\Psi_k\left(\widetilde{\mathcal{A}}^a+\widetilde{\mathcal{C}}
^a\right)\right)_ {2^j}
\Lambda^{s}\left[ \left( \psi_j^2\left(\begin{matrix}
 \pa_t u \\
  \pa_x u
\end{matrix}\right)\right)_{2^j}\right]\right\} dx
\\  = &
\int\left\{\left( \Lambda^{s}\left[\pa_a \left(
\psi_j^2\left(\begin{matrix}
 \pa_t u \\
  \pa_x u
\end{matrix}\right)\right)_{2^j}\right]\right)^T
\left(\Psi_k\left(\widetilde{\mathcal{A}}^a+\widetilde{\mathcal{C}}
^a\right)\right)_ {2^j}
\Lambda^{s}\left[ \left( \psi_j^2\left(\begin{matrix}
 \pa_t u \\
  \pa_x u
\end{matrix}\right)\right)_{2^j}\right]\right\} dx
\\  + &
\int\left\{\left( \Lambda^{s}\left[ \left(
\psi_j^2\left(\begin{matrix}
 \pa_t u \\
  \pa_x u
\end{matrix}\right)\right)_{2^j}\right]\right)^T
\pa_a\left(\Psi_k\left(\widetilde{\mathcal{A}}^a+\widetilde{\mathcal{C
}}
^a\right)\right)_{
2^j} \Lambda^{s}\left[ \left( \psi_j^2\left(\begin{matrix}
 \pa_t u \\
  \pa_x u
\end{matrix}\right)\right)_{2^j}\right]\right\} dx
\\   + &
\int\left\{\left( \Lambda^{s}\left[ \left(
\psi_j^2\left(\begin{matrix}
 \pa_t u \\
  \pa_x u
\end{matrix}\right)\right)_{2^j}\right]\right)^T
\left(\Psi_k\left(\widetilde{\mathcal{A}}^a+\widetilde{\mathcal{C}}
^a\right)\right)_ {2^j}
\Lambda^{s}\left[ \left( \pa_a\psi_j^2\left(\begin{matrix}
 \pa_t u \\
  \pa_x u
\end{matrix}\right)\right)_{2^j}\right]\right\} dx
\\ &=0.
\end{split}
\end{equation}}
Since $\widetilde{\mathcal{A}}^a$ and $\widetilde{\mathcal{C}}^a$ are
symmetric the first
and the third terms of the right hand side of (\ref{eq:4.49}) are
equal and hence {\small
\begin{equation}
\label{eq:4.50}
\begin{split}
 & \frac{2}{2^j}\mid\left\langle \Lambda^{s} \left[ \left(
\psi_j^2\left(\begin{matrix}
 \pa_t u \\
  \pa_x u
\end{matrix}\right)\right)_{2^j}\right] ,
\left(\Psi_k\left(\widetilde{\mathcal{A}}^a+\widetilde{\mathcal{C}}
^a\right)\right)_ {2^j}
\Lambda^{s}\left[ \left( \pa_a\psi_j^2\left(\begin{matrix}
 \pa_t u \\
  \pa_x u
\end{matrix}\right)\right)_{2^j}\right]\right\rangle_{L^2}\mid
\\ = &
\frac{1}{2^j}\mid \left\langle \Lambda^{s} \left[ \left(
\psi_j^2\left(\begin{matrix}
 \pa_t u \\
  \pa_x u
\end{matrix}\right)\right)_{2^j}\right] ,
\pa_a\left(\Psi_k\left(\widetilde{\mathcal{A}}^a+\widetilde{\mathcal{C
}}
^a\right)\right)_{
2^j} \Lambda^{s}\left[ \left( \psi_j^2\left(\begin{matrix}
 \pa_t u \\
  \pa_x u
\end{matrix}\right)\right)_{2^j}\right]  \right\rangle_{L^2}\mid
\\ \lesssim &
\left\lbrace
\left\|\left(
\widetilde{\mathcal{A}}^a\right)_{2^j}\right\|_{L^\infty}
+\left\|\left(
\pa_a\widetilde{\mathcal{A}}^a\right)_{2^j}\right\|_{L^\infty}
+1\right\rbrace
\left\{ \left\|(\psi_j^2\pa_t
u)_{2^j}\right\|_{H^{s}}^2+\left\|(\psi_j^2\pa_x
u)_{2^j}\right\|_{H^{s}}^2\right\}.
\end{split}
\end{equation}
That completes the estimation of the first term of the
right hand side of (\ref{eq:4.40}). The second and the third terms of
are easier to
handle since they do not contain derivatives of high order.
Recalling conditions (\ref{item:L6}) and (\ref{item:L7})
and  using algebra in $H^{s}$ for $s>\frac{3}{2}$, we have
\begin{equation}
 \label{eq:4.57}
\begin{split}
& \mid \left\langle \Lambda^{s} \left[ \left( \psi_j^2
\left( \begin{array}{c}
\pa_t u \\ \pa_x u
\end{array}\right)
\right)_{2^j}\right], \Lambda^{s}\left[ \left( \Psi_k\psi_j^2\left(
\begin{array}{cc}
{\bf b}_{22} & {\bf b}_{23} \\ {\bf b}_{32} & {\bf b}_{33}
\end{array}\right) \left( \begin{array}{c}
\pa_t u \\ \pa_x u
\end{array}\right) \right)_{2^j}\right] \right\rangle_{L^2}\mid
\\ \lesssim &
\left\|\left(
\psi_j{\widetilde{{\mathcal{B}}}}\right)_{2^j}\right\|_{H^{s}}
\left(\left\|\left( \psi_j^2
\pa_t u\right)_{2^j}\right\|_{H^{s}}+\left\|\left( \psi_j^2 \pa_x
u\right)_{2^j}\right\|_{H^{s}}\right)\\  \times &
\left(\left\|\left(\psi_j \pa_t
u\right)_{2^j}\right\|_{H^{s}}+\left\|\left(\psi_j \pa_x
u\right)_{2^j}\right\|_{H^{s}}\right)
\end{split}
\end{equation}
and
\begin{equation}
 \label{eq:4.65}
\begin{split}
& \mid \left\langle \Lambda^{s} \left[ \left( \psi_j^2
\left( \begin{array}{c}
\pa_t u \\ \pa_x u
\end{array}\right)
\right)_{2^j}\right], \Lambda^{s}\left[ \left( \Psi_k\psi_j^2\left(
\begin{array}{c}
f_2 \\ f_3
\end{array}\right)  \right)_{2^j}\right] \right\rangle_{L^2}\mid
\\ \lesssim &
\left\|\left( \psi_j^2 \pa_t u\right)_{2^j}\right\|_{H^{s}}
\left\|\left( \psi_j^2{f}_{2}\right)_{2^j}\right\|_{H^{s}}
+\left\|\left( \psi_j^2 \pa_x
u\right)_{2^j}\right\|_{H^{s}}\left\|\left(
\psi_j^2{f}_{3}\right)_{2^j}\right\|_{H^{s}}.
\end{split}
\end{equation}

To complete the proof we need to summarize
\begin{equation}
\label{eq:4.52}
\sum_{j=0}^\infty\left(\sum_{k=j-3}^{j+4}
2^{(\frac{3}{2}+\delta+1)2j}\left( E_{\pa_t}(j,k)+
E_{\pa_x}(j,k)\right)  \right).
\end{equation}
We see from the inequalities (\ref{eq:4.23}), (\ref{eq:4.24}),
(\ref{eq:4.28}), (\ref{eq:4.32}), (\ref{eq:46}), (\ref{eq:4.48}),
(\ref{eq:4.50}), (\ref{eq:4.57}) and (\ref{eq:4.65})
that the estimation of (\ref{eq:4.52})  consists  of the following
types of series:
\begin{enumerate}
\item[{\it Type 1:}]
\begin{equation}
 \label{eq:4.58}
\sum_{j=0}^\infty\left(\sum_{k=j-3}^{j+4}
2^{(\frac{3}{2}+\delta+1)2j}\left(\|(h)_{2^j}\|_{L^\infty}
\left\|(\psi_j^2f)_{2^j}\right\|_{H^{s}}
\left\|(\psi_j^2g)_{2^j}\right\|_{H^{s}}\right)  \right),
\end{equation}
where $f$ and $ g$ belong to  $ H_{s,\delta+1}$ and $h$ is in
$L^\infty$;
\item[{\it Type 2:}]
\begin{equation}
 \label{eq:4.64}
\sum_{j=0}^\infty\left(\sum_{k=j-3}^{j+4}
2^{(\frac{3}{2}+\delta+1)2j}\left(\|(h)_{2^j}\|_{L^\infty}
\left\|(\psi_j^2f)_{2^j}\right\|_{H^{s}}
\left\|(\psi_jg)_{2^j}\right\|_{H^{s}}\right)  \right),
\end{equation}
where $f$ and $ g$ belong to  $ H_{s,\delta+1}$ and $h$ is in
$L^\infty$;
\item[{\it Type 3:}]
\begin{equation}
 \label{4.55}
\sum_{j=0}^\infty\left(\sum_{k=j-3}^{j+4}
2^{(\frac{3}{2}+\delta+1)2j}\left(\|(\psi_kh)_{2^k}\|_{H^{s+1}}
\left\|(\psi_j^2f)_{2^j}\right\|_{H^{s}}
\left\|(\psi_j^2g)_{2^j}\right\|_{H^{s}}\right)  \right),
\end{equation}
where $f$ and $ g$ belong to  $ H_{s,\delta+1}$ and $h$ is in
$H_{s+1,\delta}$;
\item[{\it Type 4:}]
\begin{equation}
 \label{4.59}
\sum_{j=0}^\infty\left(\sum_{k=j-3}^{j+4}
2^{(\frac{3}{2}+\delta+1)2j}\left(\|(\psi_kh)_{2^k}\|_{H^{s+1}}
\left\|(\psi_j^2f)_{2^j}\right\|_{H^{s}}
\left\|(\psi_j g)_{2^j}\right\|_{H^{s}}\right)  \right),
\end{equation}
where $f$ and $ g$ belong to  $ H_{s,\delta+1}$ and $h$ is in
$H_{s+1,\delta}$;
\item[{\it Type 5:}]
\begin{equation}
 \label{eq:4.56}
\sum_{j=0}^\infty\left(\sum_{k=j-3}^{j+4}
2^{(\frac{3}{2}+\delta+1)2j}\left(\|(\psi_jh)_{2^j}\|_{H^{s+1}}
\left\|(\psi_j^2f)_{2^j}\right\|_{H^{s}}
\left\|(\psi_j^2g)_{2^j}\right\|_{H^{s}}\right)  \right),
\end{equation}
where $f$ and $ g$ belong to  $ H_{s,\delta+1}$ and $h$ is in
$H_{s+1,\delta}$;
\item[{\it Type 6:}]
\begin{equation}
 \label{eq:4.51}
\sum_{j=0}^\infty\left(\sum_{k=j-3}^{j+4}
2^{(\frac{3}{2}+\delta+1)2j}\left(\|(\psi_jh)_{2^j}\|_{H^{s+1}}
\left\|(\psi_j^2f)_{2^j}\right\|_{H^{s}}
\left\|(\psi_j g)_{2^j}\right\|_{H^{s}}\right)  \right),
\end{equation}
where $f\in H_{s,\delta}$, $ g\in H_{s,\delta+1}$  and $h\in
H_{s,\delta+1}$.

\end{enumerate}

The estimation (\ref{eq:4.58})- (\ref{eq:4.51}) will be done by means
of the Cauchy-Schwarz and H\"older's inequalities and the equivalence
property (\ref{eq:const:13}), of the $H_{s,\delta}$-norm. Starting
with type 1, we see from  (\ref{eq:4.58}) that

\begin{equation}
 \label{eq:4.60}
\begin{split}
 & \sum_{j=0}^\infty\sum_{k=j-3}^{j+4}
2^{(\frac{3}{2}+\delta+1)2j}\left(\left\|(\psi_j^2f)_{2^j}\right\|_{H^
{s}} \left\|(\psi_j^2g)_{2^j}\right\|_{H^{s}}\right)
\\ \leq 7 & \|h\|_{L^\infty}
\sum_{j=0}^\infty
2^{(\frac{3}{2}+\delta+1)2j}\left(\left\|(\psi_j^2f)_{2^j}\right\|_{H^
{s}}^2 + \left\|(\psi_j^2g)_{2^j}\right\|_{H^{s}}^2\right)
\\ \simeq & \|h\|_{L^\infty}\left(
\left\|f\right\|_{H_{s,\delta+1,2}}^2+
\left\|g\right\|_{H_{s,\delta+1,2}}^2\right).
\end{split}
\end{equation}
Similarly we estimate type 2, the only difference is the  use  of the
equivalence  (\ref{eq:const:13}) in the final step.
Number 3 is more sophisticated, we first note that
$(\frac{3}{2}+\delta+1)2j\leq
(\frac{3}{2}+\delta)j+(\frac{3}{2}+\delta+1)j+(\frac{3}{2}+\delta+1)j$
for $\delta\geq-\frac{3}{2}$, then we apply the H\"older inequality
with $\frac{1}{2}$, $\frac{1}{4}$ and $\frac{1}{4}$ and get
\begin{equation*}
 \label{4.62}
\begin{split} & \quad \  \sum_{j=0}^\infty\left(\sum_{k=j-3}^{j+4}
2^{(\frac{3}{2}+\delta+1)2j}\left(\|(\psi_kh)_{2^k}\|_{H^{s+1}}
\left\|(\psi_j^2f)_{2^j}\right\|_{H^{s}}
\left\|(\psi_j^2g)_{2^j}\right\|_{H^{s}}\right)  \right) \\ & \leq
\sum_{j=0}^\infty\sum_{k=j-3}^{j+4} \left(
2^{(\frac{3}{2}+\delta)j}\left\|(\psi_kh)_{2^k}\right\|_{H^{s+1}}
\right) \left(
2^{(\frac{3}{2}+\delta+1)j}\left\|(\psi_j^2f)_{2^j}\right\|_{H^{s}}
\right) \left(
2^{(\frac{3}{2}+\delta+1)j}\left\|(\psi_j^2g)_{2^j}\right\|_{H^{s}}
\right)  \\ & \leq \left( \sum_{j=0}^\infty\sum_{k=j-3}^{j+4} \left(
2^{(\frac{3}{2}+\delta)j}\left\|(\psi_kh)_{2^k}\right\|_{H^{s+1}}
\right)^2\right)^{\frac{1}{2}} \left(
\sum_{j=0}^\infty\sum_{k=j-3}^{j+4} \left(
2^{(\frac{3}{2}+\delta+1)j}\left\|(\psi_j^2f)_{2^j}\right\|_{H^{s}}
\right)^4\right)^{\frac{1}{4}}  \\ & \times \left(
\sum_{j=0}^\infty\sum_{k=j-3}^{j+4} \left(
2^{(\frac{3}{2}+\delta+1)j}\left\|(\psi_j^2g)_{2^j}\right\|_{H^{s}}
\right)^4\right)^{\frac{1}{4}} \\ & \leq
2^{(\frac{3}{2}+\delta)3}\sqrt{7}\left(
\sum_{j=0}^\infty\sum_{k=j-3}^{j+4} \left(
2^{(\frac{3}{2}+\delta)2k}\left\|(\psi_kh)_{2^k}\right\|_{H^{s+1}}
^2\right)\right)^{\frac{1}{2}} \left( \sum_{j=0}^\infty\left(
2^{(\frac{3}{2}+\delta+1)2j}\left\|(\psi_j^2f)_{2^j}\right\|_{H^{s}}
^2\right)\right)^{\frac{1}{2}}  \\ & \times \left( \sum_{j=0}^\infty
\left(
2^{(\frac{3}{2}+\delta+1)2j}\left\|(\psi_j^2g)_{2^j}\right\|_{H^{s}}^2
\right)\right)^{\frac{1}{2}}\\  & \leq
2^{(\frac{3}{2}+\delta)3}{7}^{\frac 3 2}\left( \sum_{k=0}^\infty\left(
2^{(\frac{3}{2}+\delta)2k}\left\|(\psi_kh)_{2^k}\right\|_{H^{s+1}}
^2\right)\right)^{\frac{1}{2}} \left( \sum_{j=0}^\infty\left(
2^{(\frac{3}{2}+\delta+1)2j}\left\|(\psi_j^2f)_{2^j}\right\|_{H^{s}}
^2\right)\right)^{\frac{1}{2}}  \\ & \times \left( \sum_{j=0}^\infty
\left(
2^{(\frac{3}{2}+\delta+1)2j}\left\|(\psi^2_jg)_{2^j}\right\|_{H^{s}}^2
\right)\right)^{\frac{1}{2}} \\ & \lesssim
\left\|h\right\|_{H_{s+1,\delta}}\left\|f\right\|_{H_{s,\delta+1,2}}
\left\|g\right\|_{H_{s,\delta+1,2}}  \leq
\left\|h\right\|_{H_{s+1,\delta}}\left(
\left\|f\right\|_{H_{s,\delta+1,2}}^2+\left\|g\right\|_{H_{s,\delta+1,
2}}^2\right)
\end{split}
\end{equation*}
The estimations of Types 4, 5 and 6 are similar to the last one.

We may conclude now  that
\begin{equation}
\label{eq:4.63}
\begin{split}
 & \ \ \left\langle  \left(\begin{array}{c}\pa_t u\\ \pa_x
u\end{array}\right),\pa_t\left(\begin{array}{c}\pa_t u\\ \pa_x
u\end{array}\right)\right\rangle_{s,\delta+1,{\bf a}_{33}^0}\\ &=
\sum_{j=0}^\infty\left(\sum_{k=j-3}^{j+4}
2^{(\frac{3}{2}+\delta+1)2j}\left( E_{\pa_t}(j,k)+
E_{\pa_x}(j,k)\right)  \right)
\\ & \leq
C\left(\left\|\pa_t u\right\|_{H_{s,\delta+1,2}}^2 +\left\|\pa_x
u\right\|_{H_{s,\delta+1,2}}^2 +1\right)
\leq
C\left(\| U\|_{X_{s,\delta}}^2 +1\right),
\end{split}
 \end{equation}
where the constant $C$ depends on
$\|({\mathcal{A}}^0-\e)\|_{H_{s+1,\delta}}$,
$\|{\mathcal{A}}^a\|_{H_{s+1,\delta}}$,
$\|\widetilde{{\mathcal{B}}}\|_{H_{s,\delta+1}}$,
$\|{\mathcal{F}}\|_{H_{s,\delta+1}}$,
$\|{\mathcal{A}}^\alpha\|_{L^\infty}$,
$\|\pa_x
{\mathcal{A}}^\alpha\|_{L^\infty}$, $s$ and $\delta$. By the
embeddings
Proposition \ref{prop:2.1}:\ref{Embedding} and
\ref{prop:2.1}:\ref{embedding 1}, we may replace
$\|{\mathcal{A}}^\alpha\|_{L^\infty}$ and  $\|\pa_x
{\mathcal{A}}^\alpha\|_{L^\infty}$ by their corresponding
$H_{s,\delta}$
norm.
Thus combining  inequalities (\ref{eq:4.11}), (\ref{eq:4.5}),
(\ref{eq:4.63}) with Corollary \ref{cor:Equivalence of norms}   we get
that
\begin{equation}
\frac{d}{dt}\frac{1}{2}\la
U(t),U(t)\ra_{X_{s,\delta,{\mathcal{A}}^0}}\leq
C\left(
\| U\|_{X_{s,\delta}}^2 +1\right)\leq Cc_0
\left(\| U\|_{X_{s,\delta,{\mathcal{A}}^0}}^2 +1\right),
\end{equation}
here $c_0$ is the constant of the equivalence (\ref{eq:4.2}) and in
 addition  $C$ also depends on
 $\|\pa_t {\bf
a}_{33}^0\|_{L^\infty}=\|\pa_t {\mathcal{A}}^0\|_{L^\infty}$. This
completes the proof of the energy estimates.
\hfill{$\square$}

\subsection{$L_\delta^2$ - energy estimates}
\label{subsec:energy estimates2}
The $L_\delta^2$ space is the closure of all continuous functions
with respect to the norm
\begin{equation}
 \label{eq:4.66}
\| u\|_{L_\delta^2}^2=\int(1+|x|)^{2\delta}|u(x)|^2dx.
\end{equation}
 Similarly to Definition \ref{defn:2}, we set $Y_\delta =
L^2_{\delta}\times L^2_{\delta+1}\times L^2_{\delta+1}$ and the
norm of
$V=(v_1,v_2,v_3)\in Y_\delta$ is denoted by
\begin{equation}
\label{eq:4.755}
 \|V\|_{Y_\delta}^2=\|v_1\|_{L_\delta^2}^2+\|
v_2\|_{L_{\delta+1}^2}^2+\|v_3\|_{L_{\delta+1}^2}^2.
\end{equation}
The equivalence of norms 
$\|V\|_{X_{0,\delta}}\simeq \|V\|_{Y_\delta}$
 follows from  Proposition \ref{prop:2.1}:\ref{equivalence}.

In analogous to Definition \ref{def:inner-product}, we define an
inner-product which is appropriate to the system (\ref{eq:4.34}).
So let  ${\bf a}_{33}^0$ be a positive definite matrix and $V, \Phi\in
Y_{\delta}$ be two vector valued
functions. We define an inner-product:
\begin{equation}
 \label{eq:inner-product2}
\begin{split}
\la V,\Phi\ra_{Y_{\delta},{\bf a}_{33}^0} &=\int(1+|x|)^{2\delta}
v_1^T\phi_1dx
\\ & +
\int(1+|x|)^{2\delta+2}\left[
\left(v_2^T,v_3^T\right)\left(\begin{array}{cc}\e & {\bf 0} \\ {\bf
0} &
{\bf a}_{33}^0\end{array}\right)\left(\begin{array}{ l } \phi_2 \\
\phi_3\end{array}\right)\right] dx,
\end{split}
\end{equation}
and the  norm which is associated with this product: $\|
V\|_{Y_\delta,{\bf a}_{33}^0}^2=\la
V,V\ra_{L_\delta^2,{\bf a}_{33}^0}$. If ${\bf a}_{33}^0$
satisfies (\ref{eq:4.2}), then
\begin{equation}
 \label{eq:4.2.1}
\frac 1 {c_0}\| V\|_{Y_\delta,{\bf a}_{33}^0}^2\leq
\| V\|_{Y_\delta}^2\leq c_0 \| V\|_{Y_\delta,{\bf a}_{33}^0}^2.
\end{equation}

\begin{lem}[$L_\delta^2$ energy estimates]
\label{lem:Energy eatimats2} Assume the coefficients of
(\ref{eq:4.34}) satisfy conditions (\ref{item:L1}),
(\ref{eq:4.2}), (\ref{item:L9}),  and (\ref{item:L10}). If
$U(t,\cdot)=(u(t,\cdot),\pa_t u(t,\cdot),\pa_x u(t,\cdot))\in
X_{1,\delta}$ is a solution  to the linear system
(\ref{eq:4.34}), then
\begin{equation}
 \label{eq:4.61}
\frac{d}{dt}\| U(t)\|_{Y_\delta,{\bf a}_{33}^0}^2\leq Cc_0 \left( \|
U(t)\|_{Y_\delta,{\bf a}_{33}^0}^2 +\|
\mathcal{F}\|_{L_{\delta+1}}^2\right),
\end{equation}
 where the constant $C$ depends on the $L^\infty$-norm of
${\mathcal{A}}^\alpha$,
$\pa_\alpha {\mathcal{A}}^\alpha$ and ${\mathcal{B}}$.
\end{lem}

\noindent
\textbf{Proof} (Lemma \ref{lem:Energy eatimats2}){\bf .}
Taking the derivative of $\la U(t),U(t)\ra_{L_\delta^2,{\bf
a}_{33}^0}$ yields,
\begin{equation}
 \label{eq:4.67}
\begin{split}
\frac 1 2\frac{d}{dt}\| U(t)\|_{Y_\delta,{\bf a}_{33}^0}^2 & =\la
U(t),\pa_t U(t)\ra_{Y_\delta,{\bf a}_{33}^0}
  +
\frac 1 2\int(1+|x|)^{2\delta+2}(\pa_x u)^T\pa_t {\bf a}_{33}^0(
 \pa_x u)  dx.
\end{split}
\end{equation}
By the Cauchy Schwarz inequality, the second term of the right hand
side of (\ref{eq:4.67}) is less than
\begin{equation}
 \label{eq:4.68}
\sqrt{3N}\|\pa_t {\bf a}_{33}^0\|_{L^\infty}\| \pa_x
u\|_{L_{\delta+1}}^2.
\end{equation}
Let $\widetilde{\mathcal{A}}^a$  and $\widetilde{\mathcal{C}}^a$ be
the
matrices which  is
defined  in
(\ref{eq:4.39}), since $U(t)$ satisfies system (\ref{eq:4.34}) we have
\begin{equation}
 \label{eq:4.69}
\begin{split}
&\la U(t),\pa_t U(t)\ra_{Y_\delta,{\bf a}_{33}^0}
=\int(1+|x|)^{2\delta}u^T(\pa_t u) dx
\\ & +
\sum_{a=1}^3\int(1+|x|)^{2\delta+2}
\left[ \left((\pa_t u)^T,(\pa_x
u)^T\right)\left(\widetilde{\mathcal{A}}^a+\widetilde{\mathcal{C}}
^a\right)\pa_a\left(\begin{array}{l}
\pa_t u\\ \pa_x u\end{array}\right)\right] dx
\\ & +
\int(1+|x|)^{2\delta+2}\left[\left((\pa_t
u)^T,(\pa_x u)^T\right)\left(\begin{array}{cc}{\bf b}_{22} & {\bf
b}_{23}
\\ {\bf b}_{32} & {\bf b}_{33}\end{array}\right)\left(\begin{array}{c}
 \pa_t u\\ \pa_x u\end{array}\right)\right] dx
\\ & +
\int(1+|x|)^{2\delta+2}\left[ (\pa_t u)^Tf_{2}+(\pa_x u)^T
{f}_{3}\right] dx
\\ & =:
L_1+\sum_{a=1}^3L_{2,a}+L_3+L_4.
\end{split}
\end{equation}
We are estimating each term separately:
\begin{equation}
 \label{eq:4.70}
|L_1|\leq \|u\|_{L_{\delta}^2}\|\pa_t u\|_{L_{\delta}^2}\leq \frac  1
2\left(\|u\|_{L_{\delta}^2}^2+\|\pa_t u\|_{L_{\delta+1}^2}^2\right),
\end{equation}
\begin{equation}
 \label{eq:4.73}
|2L_{2,a}|\leq 2\sqrt{4N}\left(|2\delta+2|\|
\widetilde{\mathcal{A}}^a\|_{L^\infty}+\|
\pa_a\widetilde{\mathcal{A}}^a\|_{L^\infty}+1\right)\left(\|\pa_t
u\|_{L_{\delta+1}^2}^2+\|\pa_x u\|_{L_{\delta+1}^2}^2\right),
\end{equation}
\begin{equation}
 \label{eq:4.71}
|L_3|\leq 4N\|B\|_{L^\infty}\left(\|\pa_t u\|_{L_{\delta+1}^2}^2+
\|\pa_xu\|_{L_{\delta+1}^2}^2\right)
\end{equation}
and
\begin{equation}
 \label{eq:4.72}
|L_4|\leq \frac 1 2\left(\|\pa_t
u\|_{L_{\delta+1}^2}^2+\|f_2\|_{L_{\delta+1}^2}^2+\|\pa_x
u\|_{L_{\delta+1}^2}^2+\|f_3\|_{L_{\delta+1}^2}^2\right).
\end{equation}
In (\ref{eq:4.73}) we have used the identity
\begin{equation*}
\begin{split}
0 & =\int\pa_a\left\lbrace (1+|x|)^{2\delta+2}
\left[ \left((\pa_t u)^T,(\pa_x u)^T\right)\left(
\widetilde{\mathcal{A}}^a+\widetilde{\mathcal{C}}^a\right)
\left(\begin{array}{l} \pa_t u\\
\pa_x u\end{array}\right)\right]\right\rbrace  dx
\\ & =
\int (2\delta+2) (1+|x|)^{2\delta+1}\frac{x_a}{|x|}
\left[ \left((\pa_t u)^T,(\pa_x u)^T\right)
\left(
\widetilde{\mathcal{A}}^a+\widetilde{\mathcal{C}}^a\right)
\left(\begin{array}{l} \pa_t u\\
\pa_x u\end{array}\right)\right]  dx
\\ & +
\int  (1+|x|)^{2\delta+2}
\left[ \pa_a\left(\left(\pa_t u\right)^T,\left(\pa_x u\right)^T\right)
\left(\widetilde{\mathcal{A}}^a+\widetilde{\mathcal{C}}^a\right)\left(
\begin{array}{l}
\pa_t u\\ \pa_x u\end{array}\right)\right]  dx
\\ & +
\int  (1+|x|)^{2\delta+2}
\left[ \left(\left(\pa_t u\right)^T,\left(\pa_x u\right)^T\right)
\left(\widetilde{\mathcal{A}}^a+\widetilde{\mathcal{C}}
^a\right)\pa_a\left(
\begin{array}{l}\pa_t u\\ \pa_x u\end{array}\right)\right]  dx
\\ & +
\int  (1+|x|)^{2\delta+2}
\left[ \left(\left(\pa_t u\right)^T,\left(\pa_x
u\right)^T\right)\pa_a\widetilde{\mathcal{A}}^a\left(\begin{array}{l}
\pa_t
u\\
\pa_x u\end{array}\right)\right]  dx.
\end{split}
\end{equation*}
and exploited the symmetry of $\widetilde{\mathcal{A}}^a$ and
$\widetilde{\mathcal{C}}^a$.

Summing the inequalities (\ref{eq:4.68}), (\ref{eq:4.70}), 
(\ref{eq:4.73}), (\ref{eq:4.71}) and  (\ref{eq:4.72}) and taking into
account the equivalence (\ref{eq:4.2.1} ), we get inequality
(\ref{eq:4.61}).
\hfill{$\square$}

\section{Local Existence of Quasi-linear Hyperbolic Systems}
\label{sec:Local}

Let $u:\mathbb{R}\times\mathbb{R}^3\to \mathbb{R}^N$ and  set
 $U=(u,\pa_t u,\pa_xu)$,  we consider a quasi-linear first order
hyperbolic system
\begin{equation}
 \label{eq:5.1}
{\mathcal{A}}^0(u)\pa_t U =\sum_{a=1}^3
\left({\mathcal{A}}^a(u)+{\mathcal{C}}^a\right)\pa_a
U+{\mathcal{B}}(U)U
\end{equation}
under the following conditions:

\begin{con}
\label{cond:1}  All the matrices are smooth function of their
arguments
and
\begin{enumerate}
 \item \label{item:NL1} ${\mathcal{A}}^0(u)$ , ${\mathcal{A}}^a(u)$
and
${\mathcal{C}}^\alpha$ are
symmetric matrices;
\item \label{item:NL2}   ${\mathcal{A}}^0(u)=({\bf
a}_{ij}^0(u))_{ij=1,2,3}$ is a
block matrix such that ${\bf a}_{ij}^0(u)={\bf 0}$ for $i\not =j$ and
${\bf a}_{ii}^0(u)=\e$ for $i=1,2$;
\item \label{item:NL4}   ${\mathcal{A}}^a(u)=({\bf
a}_{ij}^a(u))_{ij=1,2,3}$  are
block matrices such that ${\bf a}_{1j}^a(u)={\bf 0}$ for $j,a=1,2,3$;
\item \label{item:NL5}   ${\mathcal{C}}^a=({\bf c}_{ij}^a)_{ij=1,2,3}$
are
constant block matrices such that ${\bf c}_{1j}^a={\bf 0}$ for
$j,a=1,2,3$;
\item \label{item:NL3} ${\mathcal{B}}(U)=({\bf b}_{ij}(U))_{ij=1,2,3}$
is
a
block matrix such that ${\bf b}_{i1}(U)={\bf 0}$ and ${\bf
b}_{1,j}(U)$ are constant, $i,j=1,2,3$.
\end{enumerate}
The sizes of the blocks  are ruled  according to
(\ref{eq:4.35}).
\end{con}

Clearly the system (\ref{eq:2.3}) satisfies these assumptions. The
main result of this section is the well-posedness of the system
(\ref{eq:5.1}) in $X_{s,\delta}$-spaces.

\begin{thm}[Well-posedness of quasi-linear hyperbolic symmetric
systems]
  \label{thm:1} Let $s>\frac{3}{2}$, $\delta >
-\frac{3}{2}$,  $(f,g)\in H_{s+1,\delta}\times
H_{s,\delta+1}$ and suppose
  \begin{equation}
    \label{eq:5.23} \frac{1}{\mu}
    v^Tv\leq v^T {\bf a}_{33}^0(f)v\leq \mu v^Tv,
\qquad \forall v\in\mathbb{R}^{3N}\quad \text{and some}\ \mu\in
    \setR^+.
  \end{equation}
 Then under Assumptions \ref{cond:1} there exits a positive $T$  a
unique
  $U(t)=(u(t),\pa_t u(t),\pa_x u(t))$ a solution to
(\ref{eq:5.1}) such that $U(0,x)=(f(x),g(x),\pa_xf(x))$
 and
\begin{equation}
 \label{eq:5.24}
U\in C([0,T],X_{s,\delta}).
\end{equation}
\end{thm}

\begin{rem}
\label{rem:5.1}
 We may conclude by Mixed norm estimate \ref{mixed norm} of
Proposition \ref{prop:2.1} and (\ref{eq:5.24})  that
\begin{equation}
    \label{eq:appendix:11}
    u\in C([0,T],H_{s+1,\delta})\cap C^1([0,T],H_{s,\delta+1}).
  \end{equation}
\end{rem}

%
%

We adopt Majda's method and construct the solution through an
iteration procedure \cite{majda84:_compr_fluid_flow_system_conser}.
Similar approach was carry out in \cite{BK4}, \cite{BK5} for $s>\frac
5 2$.
Here we will examine how the special assumptions of (\ref{eq:5.1})
enable us to improve  the regularity.

\subsection{Construction of the iteration scheme}
We first note that the Embedding \ref{Embedding} of Proposition
\ref{prop:2.1} implies that the initial data $(f,g,\pa_x f)$ are
continuous,
hence there is a constant $c_0\geq 1$ and a bounded domain $G_2\subset
\mathbb{R}^N$
containing $ f$ such that
\begin{equation}
 \label{eq:5.2}
\frac 1 {c_0}v^Tv\leq v^T {\bf a}_{33}^0(u)v\leq c_0 v^Tv\quad
\text{for} \ u\in G_2.
\end{equation}
According to the density properties of $H_{s,\delta}$ (Proposition
\ref{prop:2.1}:\ref{Density}), there are sequences
$\{f^k\}_{k=0}^\infty,\{g^k\}_{k=0}^\infty\subset C_0^\infty$ and a
positive constant
$R$ such
that
%
\begin{gather}
 \label{eq:5.3}
\|(f^0,g^0,\pa_x f^0)\|_{X_{s+1,\delta}}\leq C \|(f,g,\pa_x
f) \|_{X_{s,\delta}}, \\
  \label{eq:5.5}
 \| u-f^0\|_{H_{s,\delta+1,2}}\leq {R}  \Rightarrow u\in G_2,
\end{gather}
and
\begin{equation}
 \label{eq:5.6}
\|(f^k,g^k,\pa_x f^k)-(f,g,\pa_x f)\|_{X_{s,\delta}}\leq
2^{-k}\frac{R}{{4c_0}}.
\end{equation}

The iteration scheme is defined as follows: Let
$U^0(t,x)=(f^0,g^0,\pa_x f^0)$
and  $U^{k+1}(t,x)=(u^{k+1}(t,x),\pa_t u^{k+1}(t,x),\pa_x
u^{k+1}(t,x))$ be a solution to the linear initial value problem
\begin{equation}
 \label{eq:5.7}
\left\{\begin{array}{l}
{\mathcal{A}}^0(u^k)\pa_t U^{k+1} =\sum_{a=1}^3
\left({\mathcal{A}}^a+{\mathcal{C}}^a\right)(u^{k})\pa_x
U^{k+1}+{\mathcal{B}}(U^k)U^{k+1}\\ U^{k+1}(0,x)=(f^k(x),g^k(x),\pa_x
f^k(x))\end{array}\right. .
\end{equation}
The linear theory of first order symmetric hyperbolic systems (see
e.g. \cite{JOH}) guarantees the existence of a sequence
$\{U^k(t)\}\subset C_0^\infty(\mathbb{R}^3)$. Therefore for each $k$
\begin{equation}
 \label{eq:5.8}
T_k=\sup\{T: \sup_{0<t<T}\| U^k(t)-(f^0,g^0,\pa_x
f^0)\|_{X_{s,\delta}}\leq R\}>0.
\end{equation}
  We claim that there is $T^*>0$ such that $T_k\geq T^*$ for all $k$.

\subsection{Boundedness in the $X_{s,\delta}$-norm}


\begin{lem}[Boundedness in  the norm]
 \label{lem:5.1}
There is a positive constant $T^*$  such that
\begin{equation}
\label{eq:5:12}
 \sup\{T: \sup_{0<t<T}\| U^k(t)-(f^0,g^0,\pa_x
f^0)\|_{X_{s,\delta}}\leq R\}\geq T^*\quad \text{for all}\ k.
\end{equation}
\end{lem}

\noindent
\textbf{Proof} (of Lemma \ref{lem:5.1}){\bf .}
Let $V^{k+1}=U^{k+1}-U^0$, then it
satisfies the linear
system
\begin{equation}
 \label{eq:5.9}
{\mathcal{A}}^0(u^k)\pa_t
V^{k+1}=\sum_{a=1}^3\left({\mathcal{A}}^a(u^k)+{\mathcal{C}}
^a\right)\pa_a
V^{k+1}+{\mathcal{B}}(U^k)V^{k+1}+{\mathcal{F}}^k,
\end{equation}
where \begin{equation*}{\mathcal{F}}^k=\sum_{a=1}^3
\left({\mathcal{A}}^a(u^k)+{\mathcal{C}}^a\right)\pa_a
U^0+{\mathcal{B}}(U^k)U^0
\end{equation*}
and
$V^{k+1}(0,x)=(f^{k+1}(x),g^{k+1}(x),\pa_x
f^{k+1}(x))-(f^{0}(x),g^{0}(x),\pa_x f^{0}(x))$.
At this stage we need to verify that the linear system
(\ref{eq:5.9}) meets all the requirements of the energy estimates
Lemma \ref{lem:Energy estimates}. Clearly the matrices
${\mathcal{A}}^0(u^k)$, ${\mathcal{A}}^a(u^k)$ and
${\mathcal{B}}(U^k)$
satisfy conditions (\ref{item:L1}), (\ref{eq:4.2}),
(\ref{item:L9}) and (\ref{item:L10}).

 We check now that rest of the conditions of
(\ref{eq:formulation}). From the induction hypothesis
(\ref{eq:5.8}), we have that
$\|u^{k}-f^0\|_{H_{s,\delta,2}}^2+\|\pa_x u^{k}-\pa_x
f^0\|_{H_{s,\delta+1,2}}^2\leq R^2$, therefore by  Proposition
\ref{prop:2.1}:\ref{mixed norm},
$\|u^{k}-f^0\|_{H_{s+1,\delta,2}}\leq CR$. Applying the equivalence
(\ref{eq:const:13}) and Moser type estimates (Proposition
\ref{prop:2.1}:\ref{Moser}), we have
\begin{equation*}
\begin{split}
 \|{\mathcal{A}}^0(u^k)-\e\|_{H_{s+1,\delta}} & \leq
C_1\|u^k\|_{H_{s+1,\delta,2}}\leq
C_1\left(\|u^k-f^0\|_{H_{s+1,\delta,2}}+\|f^0\|_{H_{s+1,\delta,2}}
\right)\\ & \leq C_1\left(CR+\|f^0\|_{H_{s+1,\delta,2}}\right).
\end{split}
\end{equation*}
Similarly we get $\|{\mathcal{A}}^a(u^k)\|_{H_{s+1,\delta,1}}\leq
C_2\left(CR+\|f^0\|_{H_{s+1,\delta,2}}\right)$. Here the constants
$C_1$ and $C_2$ depend on $\|u^k\|_{L^\infty}$, and
$\|{\mathcal{A}}^0-\e\|_{C^m(G_2)}$, $\|{\mathcal{A}}^a\|_{C^m(G_2)}$
respectively,  which implies that conditions (\ref{item:L4}) and
(\ref{item:L5}) hold. Having shown (\ref{item:L4}) and
(\ref{item:L5}),   we conclude from Proposition
\ref{prop:2.1}:\ref{Embedding} that ${\mathcal{A}}^\alpha(u^k)\in
C_\beta^1$
 ($\beta\geq 0$). Combing it with inequalities (\ref{eq:5.5}) and
(\ref{eq:5.8}) we get
\begin{equation*}
\begin{split}
 \|\pa_t {{\mathcal{A}}}^0(u^k)\|_{L^\infty}\leq \sup_{\bar G_2}|
\frac{\pa {\mathcal{A}}^0}{\pa u}(u)|\|\pa_t u^k\|_{L^\infty}\leq C
\|\pa_t u^k\|_{H_{s,\delta+1,2}} \\ \leq
C(\|\pa_t u^k
-g^0\|_{H_{s,\delta+1,2}}+\|g^0\|_{H_{s,\delta+1,2}})\leq
C(R+\|g^0\|_{H_{s,\delta+1,2}}),
\end{split}
\end{equation*}
this gives condition (\ref{item:L8}). In order to verify condition
(\ref{item:L6}), we denote by $\widetilde{\mathcal{B}}(U^k)$ the
non-constant
 blocks of ${\mathcal{B}}(U^k)$. Then we apply again
Moser type estimates \ref{Moser} and Algebra \ref{Algebra} of
Proposition \ref{prop:2.1},  together with induction hypothesis
(\ref{eq:5.8}) and the structure of the matrix
$\mathcal{B}$  yield
\begin{equation*}
 \|\widetilde{\mathcal{B}}(U^k)\|_{H_{s,\delta+1}}\leq C_3
\|U^k\|_{X_{s,\delta}}\leq
C_3\left(R+  \|(f^0,g^0,\pa_x f^0)\|_{X_{s,\delta}}\right).
\end{equation*}
The constant $C_3$ depends on $C^m$
norm of $\widetilde{\mathcal{B}}(U)$ taking in a bounded region of
$\mathbb{R}^{5N}$
 and $\|U^k\|_{L^\infty}$.
Finally, the $H_{s,\delta}$ estimates of ${\mathcal{A}}^\alpha(u^k)$
and
$\widetilde{\mathcal{B}}(U^k)$
with Proposition \ref{prop:2.1}:\ref{Algebra} provide an upper
bound for  $\|{\mathcal{F}}^k\|_{s,\delta+1}$. Thus we have verified
all
the conditions (\ref{eq:formulation}).


%
%

We conclude that the constant $C$ of the energy estimate
(\ref{eq:4.8}) depends only on $R$ and the initial data
$\|(f^0,g^0,\pa_x f^0)\|_{X_{s,\delta}}$, hence
\begin{equation}
\label{eq:5.10}
 \frac{d}{dt}\la
V^{k}(t),V^{k}(t)\ra_{X_{s,\delta,{\mathcal{A}}^0}}\leq
C c_0\left\{\la
V^{k}(t),V^{k}(t)\ra_{X_{s,\delta,{\mathcal{A}}^0}}+1\right\},
\end{equation}
and the constant $C$ of (\ref{eq:5.10}) is independent of $k$. By
Gronwall's inequality
\begin{equation}
\label{eq:5.11}
 \| V^{k}(t)\|_{X_{s,\delta,{\mathcal{A}}^0}}^2\leq e^{Cc_0t}\left(
\| V^{k}(0)\|_{X_{s,\delta,{\mathcal{A}}^0}}^2+Cc_0t\right).
\end{equation}
Taking into account condition 
(\ref{eq:5.6}) and
Corollary \ref{cor:Equivalence of norms}, we get from (\ref{eq:5.11})
that
\begin{equation*}
\begin{split}
&\sup_{0\leq t\leq T}\|
V^{k}(t)\|_{X_{s,\delta}}^2 \leq \\ &
e^{Cc_0T}\left\{c_0^2\left(\|(f^{k},g^{k},\pa_x
f^{k})-(f,g,\pa_x f)\|_{X_{s,\delta}}^2+ \|(f^{0},g^{0},\pa_x
f^{0})-(f,g,\pa_x f)\|_{X_{s,\delta}}^2\right)+Cc_0T\right\}
\\  &  \leq e^{Cc_0T}\left(\frac{R^2}{8}+Cc_0T\right)\leq
R^2
\end{split}
\end{equation*}
provided that $ T\leq
T^*:=\sup\{t: e^{Cc_0t}\left(\frac{R^2}{8}+Cc_0t\right)\leq R^2\}$.
\hfill{$\square$}

Having shown the boundedness of $\{U^k\}$ we may conclude by the
Compact embedding, Proposition \ref{prop:2.1}:\ref{Compact
embedding}, that
$U^k\to U$ in the $X_{s',\delta'}$-norm for any $s'<s$ and
$\delta'<\delta$. By the Mixed norm estimate \ref{mixed norm}, $u^k\to
u$ in $H_{s'+1,\delta'}$ and if we chose $\frac 3 2 <s'<s$, $-\frac
3 2<\delta'<\delta$, then the Embedding into the continuous
\ref{Embedding} implies that
\begin{equation*}
 u^k(t)\to u(t)\qquad \text{in}\ \ C^1(\mathbb{R}^3),
\end{equation*}
\begin{equation*}
 \pa_t u^k(t)\to \pa_tu(t), \ \ \pa_x u^k(t)\to \pa_xu(t) \qquad
\text{in}\ \ C(\mathbb{R}^3).
\end{equation*}
Therefore $ U(t)=(u(t), \pa_t u(t), \pa_x u(t))$ is a  solution to
system (\ref{eq:5.1}) for $0 \leq t \leq T^*$.

\subsection{Weak convergence}
\label{sec:Weak converges}

Here  we show the weak converges of $\{U^k\}$ in $X_{s,\delta}$. We
chose the simplest inner-product on $X_{s,\delta}$, that is, for 
  $V=(v_1,v_2,v_3), \Phi=(\phi_1,\phi_2,\phi_3)\in
X_{s,\delta}$, we set
\begin{equation}
\label{eq:5.12}
\begin{split}
 \la V,\Phi\ra_{X_{s,\delta}}
 & =
 \sum_{j=0}^\infty 2^{( \frac{3}{2} + \delta)2j}\left\langle
\Lambda^{s}\left(\psi_j^2 v_1\right)_{2^j},
\Lambda^{s}\left(\psi_j^2\phi_2\right)_{2^j}\right\rangle_{L^2}
\\ &+
\sum_{j=0}^\infty 2^{( \frac{3}{2} + \delta+1)2j}\left\langle
\Lambda^{s}\left(\psi_j^2 v_2\right)_{2^j},
\Lambda^{s}\left(\psi_j^2 \phi_2\right)_{2^j}\right\rangle_{L^2}
\\ & +
\sum_{j=0}^\infty 2^{( \frac{3}{2} + \delta+1)2j}\left\langle
\Lambda^{s}\left(\psi_j^2 v_3\right)_{2^j},
\Lambda^{s}\left(\psi_j^2 \phi_3\right)_{2^j}\right\rangle_{L^2}.
\end{split}
\end{equation}
This definition coincides with (\ref{eq:inner3}) in the case where
${{\mathcal{A}}}^0$ is the identity matrix.

\begin{prop}\label{prop:5.2}
Given $0<s\leq\frac{s'+s''}{2}$,
$\delta\leq \frac{\delta'+\delta''}{2}$, $V\in X_{s',\delta'}$ and $
\Phi\in
X_{s'',\delta''}$,
then

\begin{equation}
\label{eq:5.13}
 |\la V,\Phi\ra_{X_{s,\delta}}|\leq \| V\|_{X_{s',\delta'}}
 \| \Phi\|_{X_{s'',\delta''}}.
\end{equation}
\end{prop}

The proof of Proposition 5.14 appears in \cite{BK4}, \cite{BK5}
with $\delta'=\delta''=\delta$. Only a slight modification of this
proof is needed  in order to include  it to
(\ref{eq:5.13}).Therefore we leave it to the reader.

\begin{lem}[Weak Convergence]
  \label{lem:Weak_Convergence:1}
  For any $ \Phi\in X_{s,\delta}$,
  \begin{equation}
    \label{eq:Weak_Convergence:1}
    \lim_k \left\langle U^{k}(t) ,\Phi\right\rangle_{X_{s,\delta}} =
\left\langle
      U(t) ,\Phi\right\rangle_{X_{s,\delta}}
  \end{equation}
  uniformly for $0\leq t\leq T^{*}$. Consequently
  \begin{equation}
    \label{eq:Weak_Convergence:2}
    \| U(t) \|_{X_{s,\delta}} \leq \lim_k\inf
\|U^k(t)\|_{X_{s,\delta}}
  \end{equation}
  and hence the solution $U(t) $ of the initial value problem
  (\ref{eq:5.1}) belongs to $C_w\left([0,T^{*}],X_{s,\delta}\right)$,
where $ C_w$ denotes the space  of functions which are continuous in
the
weak topology.
\end{lem}

\noindent
\textbf{Proof} (of Lemma \ref{lem:Weak_Convergence:1}){\bf .}
We recall that $\| U^k(t)-U(t)\|_{H_{s',\delta'}}\to 0$ for $s'<s$
and $\delta'<\delta$.
We can pick now  $s''$  and $\delta''$ such that $s<s''$,
$s<\frac{s'+s''}{2}$,  $\delta<\delta''$   and
$\delta<\frac{\delta'+\delta''}{2}$.
Given $\Phi\in X_{s,\delta}$ and $\epsilon>0$, we may
find, by Proposition \ref{prop:2.1}:\ref{Density},
$\Phi_\epsilon\in X_{s'',\delta''}$ such that
\begin{equation}
\label{eq:5.14}
 \|\Phi-\Phi_\epsilon\|_{X_{s,\delta}}\leq \frac{\epsilon}{2R}\ \ \
\text{and}\ \ \
\|\Phi_\epsilon\|_{X_{s'',\delta''}}\leq
C(\epsilon)\|\Phi\|_{X_{s,\delta}},
\end{equation}
where  $R$ is the constant of (\ref{eq:5:12}). Writing
\begin{equation}
 \label{eq:5.15}
\begin{split}
\la U^k(t)-U(t),\Phi\ra_{X_{s,\delta}} & =\la
U^k(t)-U(t),\Phi_\epsilon\ra_{X_{s,\delta}}\\ &+\la
U^k(t)-U(t),\left(\Phi-\Phi_\epsilon\right)\ra_{X_{s,\delta}}
=:I_k+II_k,
\end{split}
\end{equation}
we have by Proposition \ref{prop:5.2} and (\ref{eq:5.14}) that
\begin{equation*}
 |I_k|\leq \|
U^k(t)-U(t)\|_{X_{s',\delta'}}C(\epsilon)\|\Phi\|_{X_{s,\delta}}\to
0.
\end{equation*}
As to  the second term of (\ref{eq:5.15}), since
$\|U^k(t)-U(t)\|_{X_{s,\delta}}\leq 2R$  by (\ref{eq:5:12}), we
get from the Cauchy-Schwarz
inequality  and (\ref{eq:5.14}) that
\begin{equation*}
 |II_k|\leq \| U^k(t)-U(t)\|_{X_{s,\delta}}\|
\Phi-\Phi_\epsilon\|_{X_{s,\delta}}\leq\frac{2R\epsilon}{2R}
=\epsilon.
\end{equation*}
Thus,
\begin{equation*}
 \limsup_{k}|\la U^k(t)-U(t),\Phi\ra_{X_{s,\delta}}|\leq\epsilon
\end{equation*}
and this completes the proof of Lemma
\ref{lem:Weak_Convergence:1}.\hfill{$\square$}

\subsection{Uniqueness}
\label{sec:uniqueness}

\begin{lem}[Uniqueness]
 \label{lem:Uniqueness}
Suppose $U(t),V(t)\in X_{s,\delta}$ are
solutions to the first order symmetric
hyperbolic system  (\ref{eq:5.1}) with initial data $(f,g)$ which
satisfy (\ref{eq:5.23}), then $U(t)\equiv V(t)$.
\end{lem}

\noindent
\textbf{Proof} (of Lemma \ref{lem:Uniqueness}){\bf .}
Put $W(t)=U(t)-V(t)$, then it satisfies the linear equation
\begin{equation}
 \label{eq:5.17}
\left\{\begin{array}{l}
{\mathcal{A}}^0(u)\pa_t W =\sum_{a=1}^3
\left({\mathcal{A}}^a(u)+{\mathcal{C}}^a\right)\pa_a
W+{\mathcal{B}}(U)W+{\mathcal{F}}\\
W(0,x)=0\end{array}\right.,
\end{equation}
where
\begin{equation}
 \label{eq:5.18}
{\mathcal{F}}=\left({\mathcal{A}}^0(u)-{\mathcal{A}}^0(v)\right)\pa_t
V+\sum_{a=1}^3
\left({\mathcal{A}}^a(u)-{\mathcal{A}}^a(v)\right)\pa_a
V+\left({\mathcal{B}}(U)-{\mathcal{B}}(V)\right)V.
\end{equation}

Since $U\in X_{s,\delta}$, ${\mathcal{A}}^\alpha(u)$, $\pa_\beta
{\mathcal{A}}^\alpha(u)$ and
${\mathcal{B}}(U)$  are bounded, we can apply  Lemma \ref{lem:Energy
eatimats2} and obtain
\begin{equation}
 \label{eq:5.19}
\frac{d}{dt}\| W(t)\|_{Y_\delta^2,{\bf a}_{33}^0(u)}^2\leq Cc_0 \left(
\|
W(t)\|_{Y_\delta^2,{\bf a}_{33}^0(u)}^2
+\|{\mathcal{F}}\|_{L_{\delta+1}^2}^2\right)
\end{equation}
(See (\ref{eq:inner-product2}) for the definition of the norm $\|
W\|_{Y_\delta^2,{\bf a}_{33}^0(u)}^2$).

We turn now  to the estimation of 
$\|{\mathcal{F}}\|_{L_{\delta+1}}^2$
in terms of the difference $\|U-V\|_{Y,\delta,{\bf a}_{33}^0(u)}$. 
From the
structure of the matrices ${\mathcal{A}}^\alpha(u)$ in Assumptions
\ref{cond:1}, we
see that $\left( {\mathcal{A}}^0(u)-{\mathcal{A}}^0(v)\right) \pa_t
V=\left({\bf
a}_{33}^0(u)-{\bf a}_{33}^0(v)\right)\pa_t \pa_x v$ and
\begin{equation*}
\left(
{\mathcal{A}}^a(u)-{\mathcal{A}}^a(v)\right) \pa_a
V=\left(\widetilde{\mathcal{A}}^a(u)-\widetilde{\mathcal{A}}
^a(v)\right)\pa_a\left(
\begin{array}{c} \pa_t v\\ \pa_x v\end{array}\right),
 \end{equation*}
where $$\widetilde{\mathcal{A}}^a(p)=\left( \begin{array}{cc}{\bf
a}_{22}^a(p) &{\bf a}_{23}^a(p)\\ {\bf a}_{32}^a(p)&
{\bf a}_{33}^a(p)
\end{array}\right), \qquad a=1,2,3. $$
Our idea is to use inequality (\ref{eq:elliptic part:12})  with
$s=0$, $s_1=1$ and $s_2=s-1$   and then to apply the Difference
estimate
\ref{Difference estimate} of Proposition \ref{prop:2.1}. We
also note that by Proposition \ref{prop:2.1}:\ref{equivalence},
$\|u-v\|_{H_{1,\delta}}^2\simeq
\|u-v\|_{L^2_{\delta}}^2+\|\pa_x(u-v)\|_{L^2_{\delta+1}}^2\leq
\|U-V\|_{Y_\delta}^2$. These yield the following estimations:

\begin{equation}
 \label{eq:5.20}
\begin{split}
& \left\|\left({\bf
a}_{33}^0(u)-{\bf a}_{33}^0(v)\right)\pa_t \pa_x
v\right\|_{L_{\delta+1}^2}^2   \simeq  \left\|\left({\bf
a}_{33}^0(u)-{\bf a}_{33}^0(v)\right)\pa_t \pa_x
v\right\|_{H_{0,\delta+1}}^2 \\& \leq C   \left\| \left({\bf
a}_{33}^0(u)-{\bf a}_{33}^0(v)\right)\right\|_{H_{1,\delta}}^2 \left\|
\pa_t \pa_x v\right\|_{H_{s-1,\delta+2}}^2 \\
 & \leq   C  (\|u\|_{H_{s+1,\delta}},\|v\|_{H_{s+1,\delta}})  \left\|
u-v\right\|_{H_{1,\delta}}^2 \left\| \pa_t
v\right\|_{H_{s,\delta+1}}^2   \\ & \leq C
(\|u\|_{H_{s+1,\delta}},\|v\|_{H_{s+1,\delta}})\left\|
V \right\|_{X_{s,\delta}}^2\left\| U-V\right\|^2_{Y_\delta}.
\end{split}
\end{equation}
Similarly,
\begin{equation}
 \label{eq:5.21}
\begin{split}
&\left\|\left(\widetilde{\mathcal{A}}^a(u)-\widetilde{\mathcal{A}}
^a(v)\right)\pa_a\left(
\begin{array}{c}\pa_t v\\
\pa_x v\end{array}\right)\right\|_{L_{\delta+1}^2}
\\  \leq  C &  (\|u\|_{H_{s+1,\delta}},\|v\|_{H_{s+1,\delta}})
\|u-v\|_{H_{1,\delta}}^2
\left(\left\|  \pa_t v\right\|_{H_{s,\delta+1}}^2+\left\|  \pa_x
v\right\|_{H_{s,\delta+1}}^2\right) \\ \leq  C &
(\|u\|_{H_{s+1,\delta}},\|v\|_{H_{s+1,\delta}})
\left\|  V\right\|_{X_{s,\delta}}^2 \|U-V\|_{Y_{\delta}}^2 .
\end{split}
\end{equation}
Writing  ${\mathcal{B}}={\mathcal{B}}(p,q,r)$, then by Assumptions
\ref{cond:1}:\ref{item:NL3}
we have
\begin{equation*}
 \left(
{\mathcal{B}}(U)-{\mathcal{B}}(V)\right)V=\pa_t(u-v)\cdot\nabla_q{
\mathcal {B}}
V+\pa_x(u-v)\cdot\nabla_r{\mathcal{B}}V.
\end{equation*}
Hence  the simple  weighted $L^2$ estimate gives
\begin{equation}
 \label{eq:5.26}
\begin{split}
\left\|\pa_t(u-v)\nabla_q{\mathcal{B}}V\right\|_{L^2_{\delta+1}}^2 &
\leq
\left\|\nabla_q{\mathcal{B}}\right\|_{L^\infty}
^2\left\|\pa_t(u-v)\right\|_ { L^2_
{\delta+1}}
^2 \left\|V\right\|_{L^\infty}^2
\end{split}
\end{equation}
and
\begin{equation}
 \label{eq:5.27}
\begin{split}
\left\|\pa_x(u-v)\nabla_r{\mathcal{B}}V\right\|_{L^2_{\delta+1}}^2
\leq
\left\|\nabla_r{\mathcal{B}}\right\|_{L^\infty}
^2\left\|\pa_x(u-v)\right\|_ { L^2_
{\delta+1} }
^2 \left\|V\right\|_{L^\infty}^2.
\end{split}
\end{equation}
Thus, inequalities (\ref{eq:5.20})-(\ref{eq:5.27}) with the
equivalence
(\ref{eq:4.2.1})  show that
\begin{equation}
 \label{eq5:18}
\left\| {\mathcal{F}}\right\|_{L^2_{\delta+1}}\leq
C\left\|V
\right\|_{X_{s,\delta}}\left\|U-V\right\|_{Y_\delta}\leq C\left\|V
\right\|_{X_{s,\delta}}\left\|U-V\right\|_{Y_\delta,{\bf
a}_{33}^0(u)}.
\end{equation}

Inserting (\ref{eq5:18})  in (\ref{eq:5.19}) and using
Gronwall's
inequality we get that
\begin{equation*}
 \left\|W(t)\right\|_{Y_{\delta,{\bf a}_{33}^0(u)}}^2\leq
e^{Cc_0t}\left\|W(0)\right\|_{Y_{\delta,{\bf a}_{33}^0(u)}}^2
\end{equation*}
and since $W(0)=0$, it implies that $W(t)\equiv 0$. \hfill{$\square$}

\subsection{Continuation in the norm}
\label{sec:Continuation in the norm}

\begin{lem}[Continuation in the norm]
 \label{lem:Continuation in the norm}
Let $U(t)$ be a solutions to the first order symmetric
hyperbolic system  (\ref{eq:5.1}) with initial data $(f,g)$ which
satisfy (\ref{eq:5.23}), then (\ref{eq:5.24}) holds.
\end{lem}

\noindent
\textbf{Proof} (of Lemma \ref{lem:Continuation in the norm}){\bf .}
Since $X_{s,\delta}$ is a Hilbert space it suffices to show
that
\begin{equation*}
 \limsup_{t\to
0^+}\left\|U(t)\right\|_{X_{s,\delta,{\mathcal{A}}^0(f)}}
\leq\left\|U(0)\right\|_ { X_ { s ,\delta,{\mathcal{A}}^0(f)}}.
\end{equation*} 
Having proved the uniqueness, we may assume that $U$ is
the limit of the iteration sequence $U^k$. Furthermore, since
$u^k(t)\to u(t)$ uniformly in $[0,T^*]$ and the matrix
${\mathcal{A}}^0$
depends
solely on $u$, we see from the inner product (\ref{eq:inner3}) that
for a given $\epsilon>0$ there is a positive integer $k_0$ such that
\begin{equation}
 \label{eq:5.30}
\left\|V\right\|_{X_{s,\delta,{\mathcal{A}}^0(u(t))}}\leq
(1+\epsilon)\left\|V\right\|_{X_{s,\delta,{\mathcal{A}}^0(u^k(t))}},
\qquad k\geq
k_0,\ \ V\in X_{s,\delta}.
\end{equation}

Using the fact that  $u(t,\cdot)\to f(\cdot)$  uniformly as $t\to0$,
Lemmas (\ref{lem:Weak_Convergence:1}) and (\ref{lem:Energy
estimates}), and (\ref{eq:5.6}) we get
\begin{equation*}
 \begin{split}
  \limsup_{t\to 0^+} \left\|
U(t)\right\|_{X_{s,\delta,{\mathcal{A}}^0(f)}}^2 &=
\limsup_{t\to 0^+} \left\|
U(t)\right\|_{X_{s,\delta,{\mathcal{A}}^0(u(t))}}^2\\&\leq
\limsup_{t\to 0^+}\left(\liminf_k \left\|
U^{k+1}(t)\right\|_{X_{s,\delta,{\mathcal{A}}^0(u(t))}}^2\right)\\
&\leq
(1+\epsilon)^2
\limsup_{t\to 0^+}\left(\liminf_k \left\|
U^{k+1}(t)\right\|_{X_{s,\delta,{\mathcal{A}}^0(u^k(t))}}^2\right)\\
&\leq
(1+\epsilon)^2
\limsup_{t\to 0^+}\left(\liminf_k   e^{Cc_0t}\left(\left\|
U^{k+1}(0)\right\|_{X_{s,\delta,{\mathcal{A}}^0(u^k(0))}}
^2+Cc_0t\right)\right)\\ &
\leq (1+\epsilon)^2 \left\|
U(0)\right\|_{X_{s,\delta,{\mathcal{A}}^0(f)}}^2.
 \end{split}
\end{equation*}
This completes the proof of the Lemma and thereby of Theorem
\ref{thm:1}.   \hfill{$\square$}

\section{Proof of the main result }
\label{sec:main}

The solution of the constraint equations (\ref{eq:3}) in the weighted
Sobolev spaces of fractional order $H_{s,\delta}$ has been proved by
Maxwell
\cite{maxwell06:_rough_einst} for $s>\frac 1 2$ and Brauer and Karp
for $\geq 1$  \cite{BK5} (see also \cite{ICH12}).
Thus for a given set of free data $\left( \bar\h_{ab},\bar\K_{ab}
\right)$ such
that $\left( \bar\h_{ab}-\e_{ab},\bar\K_{ab} \right)\in
H_{s+1.\delta}\times
H_{s,\delta+1}$, there is conformally equivalent data $({\h}_{ab},
\K_{ab})$
which satisfies the constraint equations (\ref{eq:3}). Moreover, there
is a constant $C$ such that
\begin{equation}
\label{eq:6.1}
 \left\|\left( \h_{ab}-\e_{ab},\K_{ab}
\right)\right\|_{H_{s+1,\delta}\times
H_{s,\delta+1}}\leq C \left\|\left(
\bar{\h}_{ab}-\e_{ab},\bar{\K}_{ab}
\right)\right\|_{H_{s+1,\delta}\times
H_{s,\delta+1}}.
\end{equation}

We apply now Theorem \ref{thm:1} to $\left(\g_{\alpha\beta}-{\bf
m}_{\alpha\beta},\pa_t \g_{\alpha\beta},\pa_x\g_{\alpha\beta}\right)$
with initial data (\ref{eq:2}) and
where the pair $\left( \h_{ab}, \K_{ab} \right)$  satisfies the
constraint equations (\ref{eq:3}). Then $\g_{\alpha\beta}(t)$ is the
unique solution to the reduced Einstein equation (\ref{wave
equation}) and (\ref{eq:13}) holds by Remark \ref{rem:5.1}.
Inequality (\ref{eq:15}) follows from (\ref{eq:6.1}) since  for
$t\in[0,T]$ the bounds
of $\| \g_{\alpha\beta}(t)-{\bf m}_{\alpha\beta} \|_{H_{s+1,\delta}}$
and
$\|\pa_t\g_{\alpha\beta}(t)\|_{H_{s,\delta+1}}$  depend solely on the
initial data $\left( \h_{ab}, \K_{ab} \right)$.

In order to assure that $\g_{\alpha\beta}(t)$ satisfies the vacuum
Einstein equation (\ref{eq:1}) we need to establish the harmonic
condition (\ref{eq:9}). Recalling that  $F^\mu$ satisfies  the linear
wave equation
\begin{equation}
\label{}
 \g^{\alpha\beta}\pa_\alpha\pa_\beta
F^\mu-\Gamma_{\alpha\beta}^\nu\g^{\alpha\beta}\pa_\nu F^\nu=0
\end{equation}
(see e.g. \cite{CHY}, \cite{wald84:_gener_relat}), it thus suffices
to show that $F^\mu(0,x)=\pa_t F^\mu(0,x)=0$. Hence by the uniqueness
of linear hyperbolic systems, it follows that $F^\mu\equiv 0$. Note
that $\g^{\alpha\beta}-\e^{\alpha\beta}\in H_{s+1,\delta}$,
$\Gamma_{\alpha\beta}^\nu\g^{\alpha\beta}\in H_{s,\delta+1}$ and
$s>\frac 3 2$, therefore these facts  allow us to use known
uniqueness results for
linear hyperbolic symmetric systems with   coefficients in $H^s$
\cite{FMA}, \cite{Kato_70}, or alternatively, we apply the
$L_\delta^2$-energy estimate Lemma \ref{lem:Energy eatimats2},
combined with
Gronwall's inequality. We can now use the free data
$\pa_t\g_{0\alpha}$ to get the condition $F^\mu(0,x)=0$. Then
exploiting
the fact that $\left( \h_{ab}, \K_{ab} \right)$ satisfies the
constraint equations (\ref{eq:3}) leads to the second condition $\pa_t
F^\mu(0,x)=0$, see e.g. \cite{bartnik04}, \cite{wald84:_gener_relat}.

\vskip 5mm
\noindent
\textbf{Acknowledgement:} \noindent
I wish to acknowledge Uwe Brauer for many long and helpful
discussions and to Greg Galloway who enlightened the uniqueness  issue
for me.
Part of this project was done during the author's visit at the
Department
of Mathematics at Potsdam Universit\"at and I would like to thank
Professor B.-W.  Schulze for his support and kind hospitality.

\bibliographystyle{plain}
\bibliography{bibgraf}

\noindent \textsc{Department of Mathematics\\ ORT Braude College,
P.O. Box 78\\ 21982 Karmiel, Israel}\\
\textit{E-mail:}  \textsf{karp@braude.ac.il}
\\ \textit{URL:} \textsf{http://brd4.braude.ac.il/~karp/}

\end{document}